\documentclass[conference]{./IEEEtran}
\usepackage{eso-pic}

\makeatletter

\usepackage{amsmath}
\usepackage{graphicx} 
\usepackage{caption}
\usepackage{subcaption}
\usepackage{booktabs} 
\usepackage{array} 
\usepackage{paralist} 
\usepackage{verbatim} 
\usepackage{alltt}
\usepackage{subfig} 
\usepackage{dsfont}
\usepackage{mathtools}
\usepackage{cleveref}
\usepackage{cite}
\usepackage{amssymb}

\newtheorem{theorem}{Theorem}

\newtheorem{lemma}{Lemma}

\newcommand{\enma}[1]   {\ensuremath{#1}}

\newcommand{\edoc}{\end{document}}
\newcommand{\beq}{\begin{equation}}  
\newcommand{\eeq}{\end{equation}}    
\newcommand{\bea}{\begin{align}}
\newcommand{\eea}{\end{align}}
\newcommand{\beqn}{\begin{eqnarray}}
\newcommand{\eeqn}{\end{eqnarray}}
\newcommand{\beqns}{\begin{eqnarray*}}
\newcommand{\eeqns}{\end{eqnarray*}}
\newcommand{\bct}{\begin{center}}
\newcommand{\ect}{\end{center}}
\newcommand{\btmz}{\begin{itemize}}
\newcommand{\etmz}{\end{itemize}}
\newcommand{\benum}{\begin{enumerate}}
\newcommand{\eenum}{\end{enumerate}}

\newcommand{\1}{\mathbf{1}}
\newcommand{\0}{\mathbf{0}}
\newcommand{\kp}{\kappa}

\newcommand{\Norm}[1]{\left\| #1 \right\|}

\newcommand{\blockdiag} {\enma{\mathrm{blockdiag}}}

\newcommand{\image}{\mbox{\rm{Im\kern.043333em}}}

\title{Skewless Network Clock Synchronization}

\author{\IEEEauthorblockN{Enrique Mallada$^*$, Xiaoqiao Meng$^\dagger$, Michel Hack$^\dagger$, Li Zhang$^\dagger$, and Ao Tang$^*$}
\IEEEauthorblockA{$^*$ School of ECE, Cornell University, Ithaca, NY 14853, USA.
\\
$^\dagger$ IBM T. J. Watson Research Center. 1101 Kitchawan Rd, Yorktown Heights, NY 10598, USA.
\\ Email: \{em464@,atang@ece.\}cornell.edu, \{xmeng,hack,zhangli\}@us.ibm.com
}
}

\begin{document}

\maketitle

\begin{abstract}
This paper examines synchronization of computer clocks connected via a data network and proposes a skewless algorithm to synchronize them. Unlike existing solutions, which either estimate and compensate the frequency difference (skew) among clocks or introduce offset corrections that can generate jitter and possibly even backward jumps, our algorithm achieves synchronization without these problems. We first analyze the convergence property of the algorithm and provide necessary and sufficient conditions on the parameters to guarantee synchronization. We then implement our solution on a cluster of IBM BladeCenter servers running Linux and study its performance. In particular, both analytically and experimentally, we show that our algorithm can converge in the presence of timing loops. This marks a clear contrast with current standards such as NTP and PTP, where timing loops are specifically avoided. Furthermore, timing loops can even be beneficial in our scheme. For example, it is demonstrated that highly connected subnetworks can collectively outperform individual clients when the time source has large jitter. It is also experimentally demonstrated that our algorithm outperforms other well-established software-based solutions such as the NTPv4 and IBM Coordinated Cluster Time (IBM CCT).
\end{abstract}

%
\IEEEpeerreviewmaketitle

\section{Introduction}\label{sec:introduction}

Keeping consistent time among different nodes in a network is a fundamental requirement of many distributed applications. Their internal clocks are usually not accurate enough and tend to drift  apart from each other over time, generating inconsistent time values. Network clock synchronization allows these devices to correct their clocks to match a global reference of time, such as the Universal Coordinated Time (UTC), by performing  time measurements through the network. For example, for the Internet, network clock synchronization has been an important subject of research and several different protocols have been proposed~\cite{mills_network_2006, sobeih_almost_2007, _ieee_2008, veitch_robust_2009, carli_networked_2010, mallada_distributed_2011, ridoux_case_2012}.
These protocols are used for various legacy and emerging applications with diverse precision requirements such as banking
transactions, communications, traffic measurement and security protection. In particular, in modern wireless cellular networks, 
 time-sharing protocols need an accuracy of several microseconds to guarantee the efficient 
use of channel capacity. Another example is the recently announced Google Spanner~\cite{corbett_spanner_2012}, a globally-distributed database, which depends on globally-synchronized clocks within at most several milliseconds drifts. 

There are two major difficulties that make the network clock synchronization problem challenging.  First, the frequency of hardware clocks is sensitive to temperature and is constantly varying. Second, the latency introduced by the OS and network congestion delay results in errors in the time measurements. Thus, most protocols introduce different ways of estimating the frequency mismatch (skew)\cite{zhang_clock_2002, kim_tracking_2012} and measuring the time difference (offset) \cite{elson_fine-grained_2002,hunt_network_2010}.
This leads to extensive literature on skew estimation \cite{kim_tracking_2012,marouani_internal_2008,moon_estimation_1999,lemmon_model-based_2000}  which suggests that explicit skew estimation is necessary for clock synchronization.


This paper takes a different approach and shows that focusing on skew estimation could be misleading. We provide a simple algorithm that is able to compensate the clock skew without any explicit estimation of it. Our algorithm only uses current offset information and an exponential average of the past offsets. Thus, it neither needs to store long offset history nor perform expensive computations on them.
We analyze the convergence property of the algorithm and provide necessary and sufficient conditions for synchronization. The parameter values that guarantee synchronization depend on the network topology, but there exists a subset of them that is independent of topology and therefore of great practical interest.

We also discover a rather surprising fact. A common practice in the clock synchronization community is to avoid timing loops in the network~ \cite[p. 3]{mills_network_2006} \cite[p. 16, s. 6.2]{_ieee_2008}.  This is because it is thought that timing loops can introduce instability as stated in \cite{mills_network_2006}:
{\it
"Drawing from the experience of the telephone industry, which learned such lessons at considerable cost, the subnet topology...
 must never be allowed to form a loop."
} Even though for some parameter values loops can produce instability, we show that a set of proper parameters can guarantee convergence even in the presence of loops. Furthermore, we experimentally demonstrate in Section \ref{sec:experimental_results} that timing loops among clients can actually help reduce the jitter of the synchronization error and is therefore desirable.

\subsection{Related Work and Contribution}\label{sec:related_work}

Clock synchronization on computer networks has been subject of study for more than 20 years. The current {\it de facto} standard for IP networks
is the Network Time Protocol (NTP) proposed by David Mills~\cite{mills_network_2006}. It is a low-cost, purely software-based
solution whose accuracy mostly ranges from hundreds of microseconds to several milliseconds, which is often not sufficient. On the other hand,   
IEEE 1588 (PTP)~\cite{_ieee_2008}  gives superior performance by achieving  sub-microsecond or even nanosecond accuracy. However, it is relatively expensive as it requires special hardware support to achieve those accuracy levels and may not be fully compatible with legacy cluster systems. 

More recently, new synchronization protocols have been proposed with the objective of balancing between accuracy and cost. For example,
 IBM Coordinated Cluster Time (CCT)~\cite{froehlich_achieving_2008} is able to provide better performance than NTP without additional hardware. Its success is based on a skew estimation mechanism~\cite{zhang_clock_2002} that progressively adapts the clock frequency without offset corrections.
Another alternative is the RADclock~\cite{veitch_robust_2009,ridoux_case_2012} which estimates the skew and produces offset corrections, but provides a secondary relative clock that is more robust to jitter.
 
 The solution provided in this paper solves problems present on IBM CCT and RADclock. We are able to achieve microsecond level accuracy without requiring any special hardware as the previous solutions.  However, our protocol does not explicitly estimate the skew, which makes the implementation simpler and more robust to jitter than IBM CCT, and does not introduce offset corrections, which avoids the need of a secondary clock as in RADclock.  Furthermore, we present a theoretical analysis of its behavior in network environments that unveils some rather surprising facts.

The rest of the paper is organized as follows. In Section \ref{sec:comp_clocks} we provide some background on how clocks are actually implemented in computers and how different protocols discipline them. Section \ref{sec:algorithm} motivates and describes our algorithm together with an intuitive explanation of why it works. In Section \ref{sec:analysis}, we analyze the algorithm and determine the set of parameter values and connectivity patterns under which synchronization is guaranteed. Experimental results evaluating the performance of the algorithm are presented in Section \ref{sec:experimental_results}. We conclude in Section \ref{sec:conclusions}.

\section{Synchronization of Computer Clocks} \label{sec:comp_clocks}

Most computer architectures keep their own estimate of time using a counter that is periodically increased by either hardware or kernel's interrupt service routines (ISRs). On Linux platforms for instance, there are usually several different clock devices that can be selected as the clock source by changing the $clocksource$ kernel parameter. 
One particular counter that has recently been used by several clock synchronization protocols~\cite{froehlich_achieving_2008,veitch_robust_2009} is the Time Stamp Counter (TSC) that counts the number of CPU cycles since the last restart of the system. For example, in the IBM BladeCenter LS21 servers, the TSC is a $64$-bit counter that increments every $\delta^o = 0.416$ns since the CPU nominal frequency $f^o=1/\delta^o=2399.711$MHz.

\begin{figure}[htp]
       \begin{subfigure}[b]{0.49\columnwidth}
               \centering
               \includegraphics[width=.95\columnwidth,height=.65\columnwidth]{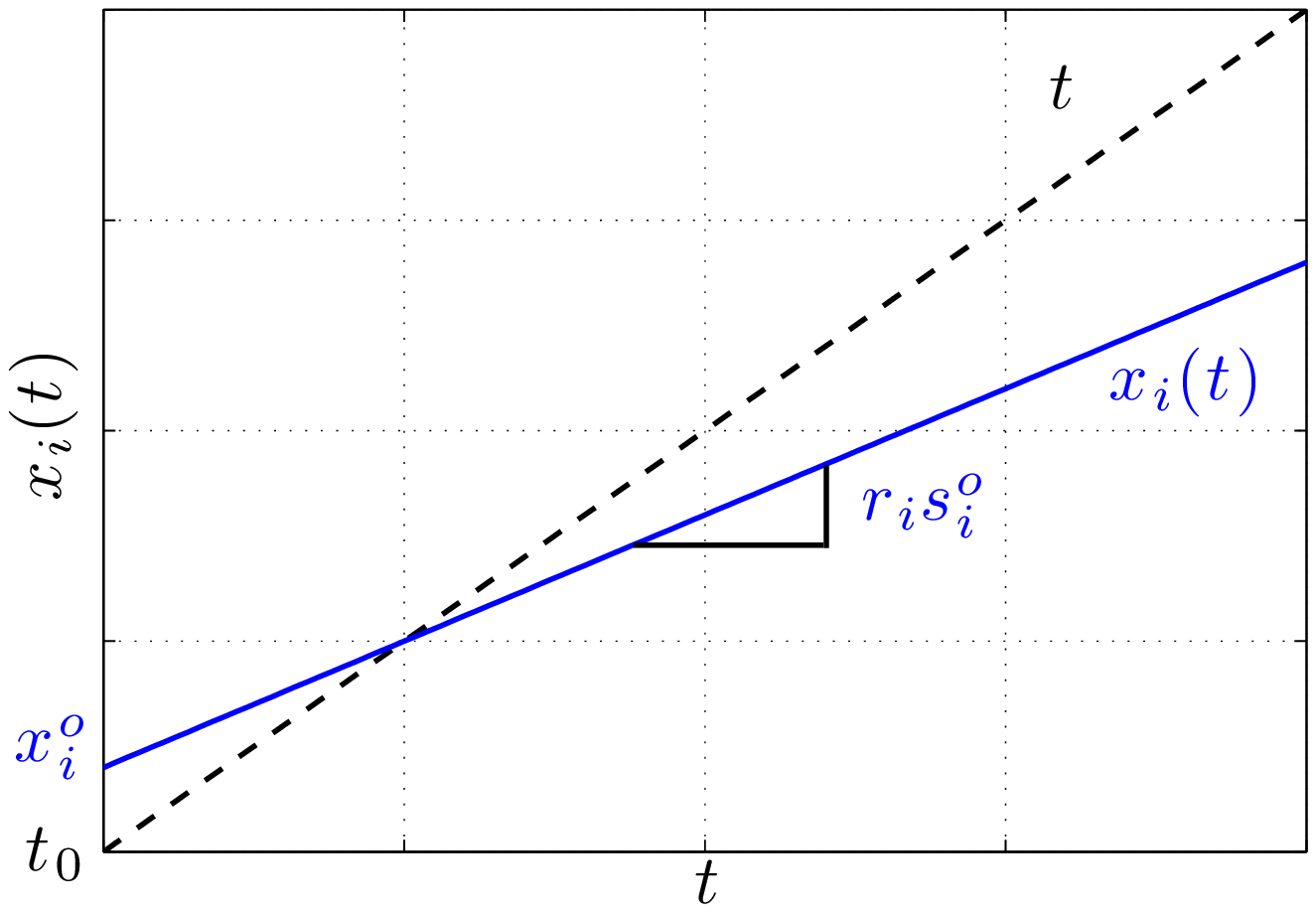}
               \caption{Illustration of computer time estimate $x_i(t)$ and UTC time $t$}\label{fig:xt}
        \end{subfigure}
        \begin{subfigure}[b]{0.49\columnwidth}
               \centering
               \includegraphics[width=.95\columnwidth,height=.65\columnwidth]{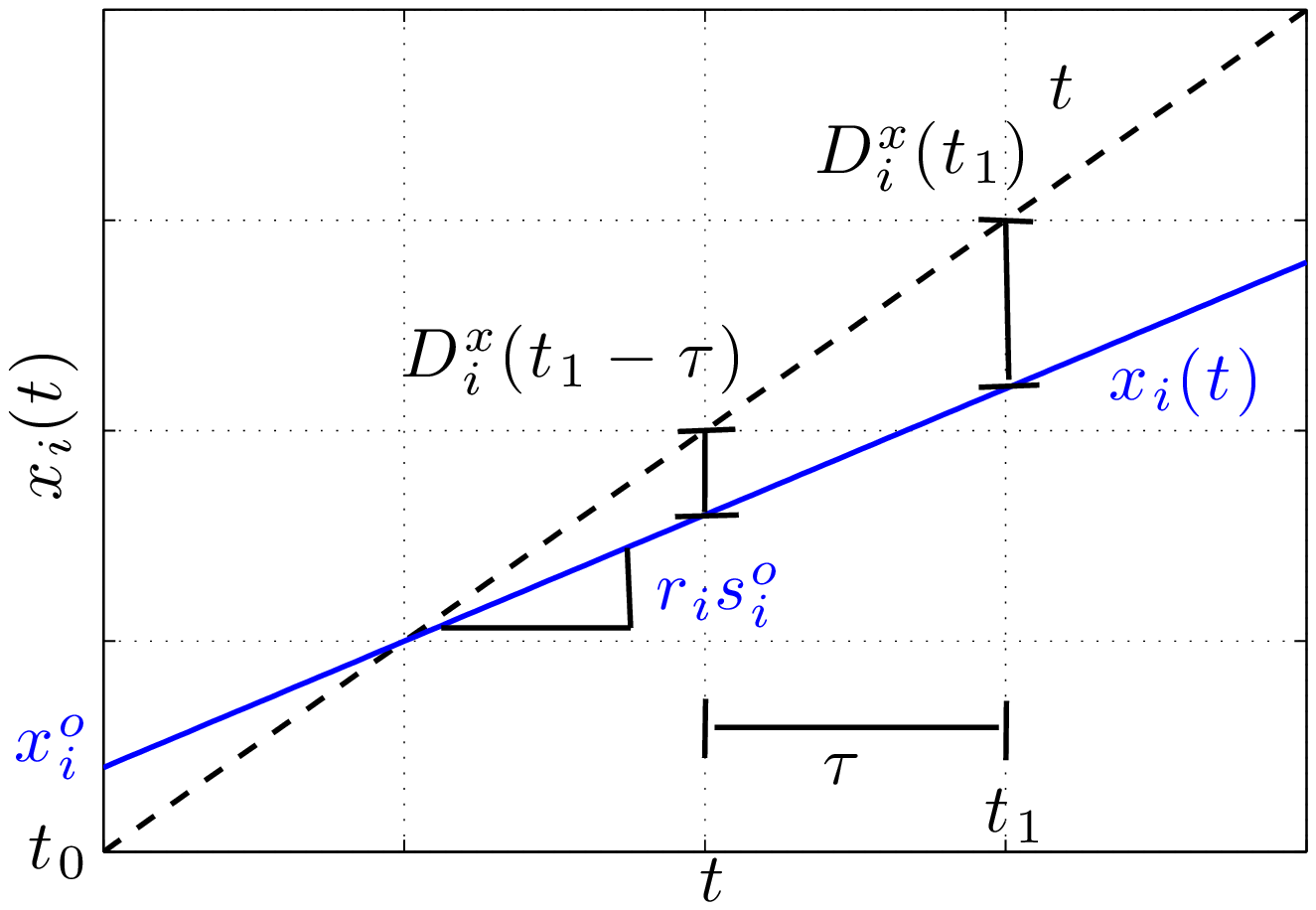}
		\caption{ Offset and relative skew measurements }\label{fig:xt3}
        \end{subfigure}%
        \vspace{-.25cm}
        \caption{Time estimation and relative measurements}\label{fig:fig1}
\end{figure}

Based on this counter, each server builds its own estimate $x_i(t)$ of the global time reference, UTC, denoted here by $t$. Thus, synchronizing computer clocks implies correcting $x_i(t)$ in order to match $t$, i.e. enforcing $x_i(t)=t$.
There are two difficulties on this estimation process. Firstly, the initial time $t_0$ in which the counter starts its unknown. Secondly,   
the clock frequency is usually unknown with enough precision and therefore presents a skew  $r_i=\frac{x_i(t)-x_i(t_0)}{t-t_0}$.
This is illustrated in Figure \ref{fig:xt} where $x_i(t)$ not only increases at a different rate than $t$, but also starts from a value different from $t_0$, represented by  $x_i^o$. 

Mathematically, $x_i(t)$ can be described by the linear map of the global reference $t$, i.e.
\begin{equation}\label{eq:linear_map}
x_i(t) = r_is_i^o(t-t_0) +x_i^o,
\end{equation}
where $s_i^o$ is an additional skew correction implemented to compensate the skew $r_i$; in Figure \ref{fig:xt} $s_i^o=1$.  Equation \eqref{eq:linear_map} also shows that if one can set $s_i^o=1/r_i$ and $x_i^o=t_0$, then we obtain a perfectly synchronized clock with $x_i(t)=t$.


\begin{figure}[htp]
       \begin{subfigure}[b]{0.49\columnwidth}
               \centering
               \includegraphics[width=\columnwidth,height = .6\columnwidth]{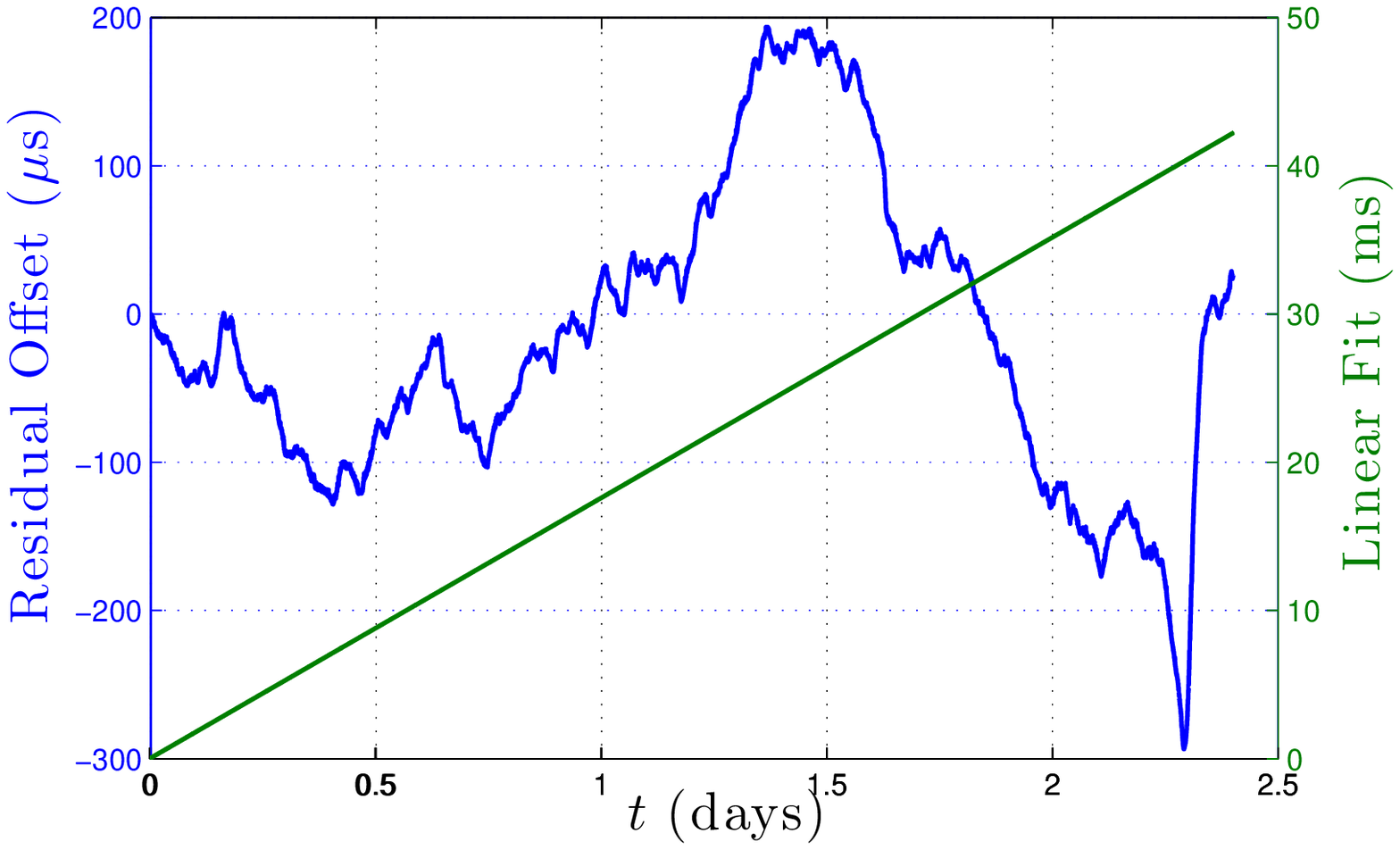}
               \caption{Offset between two TSC counters: The straight line is a linear fit that is subtracted from the offset values in residual offset axis}\label{fig:tsc_drift}
        \end{subfigure}
        \begin{subfigure}[b]{0.49\columnwidth}
               \centering
		\includegraphics[width=\columnwidth,height = .6\columnwidth]{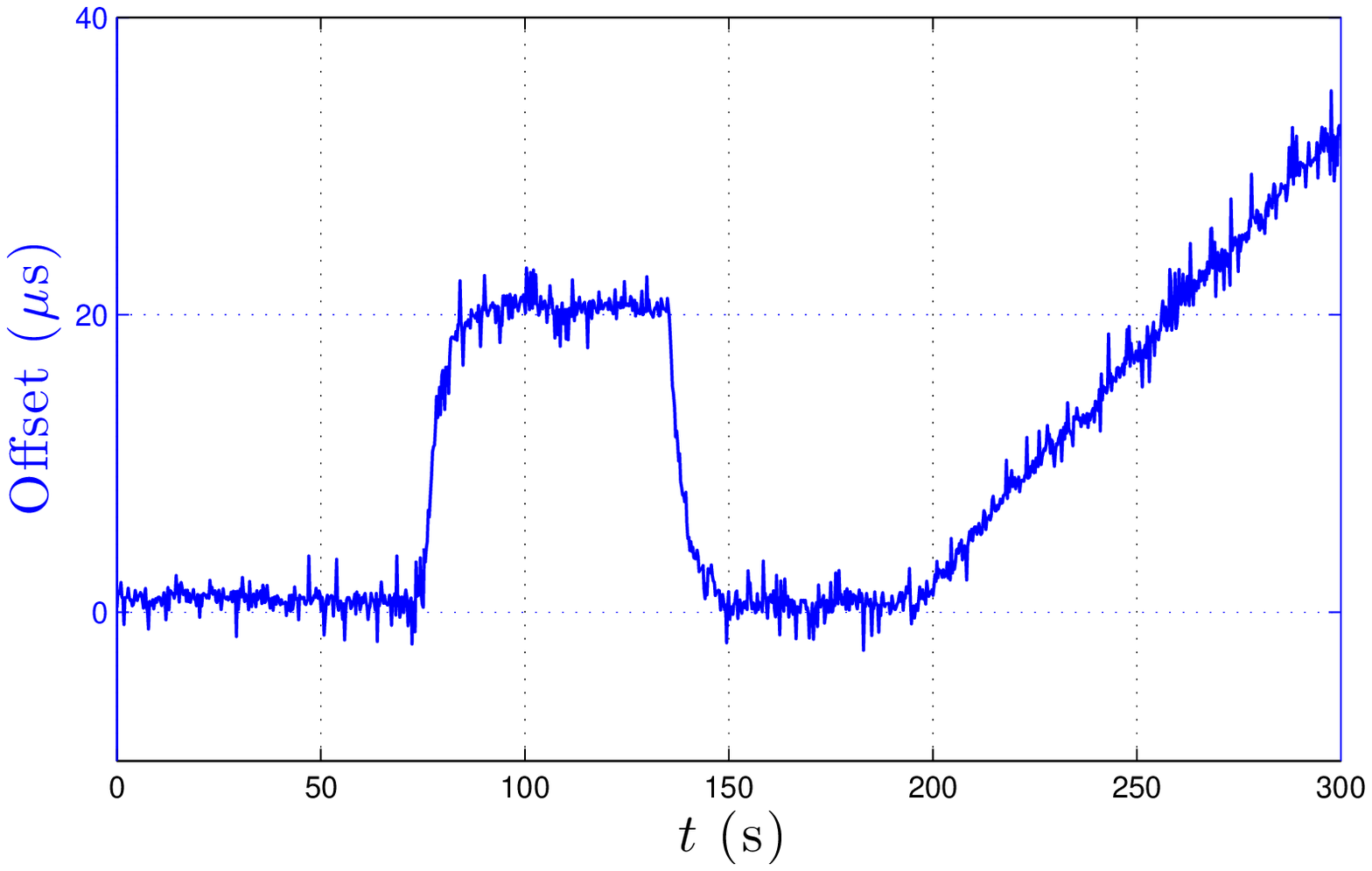}
		\caption{ Example of skew and offset corrections on linux time: First a $20\mu$s offset is added and subtracted and then a skew of $0.3$ppm is introduced}\label{fig:adjtimex}
        \end{subfigure}%
        \vspace{-.25cm}
        \caption{Comparison between two TSC counters, and skew and offset corrections using {\it adjtimex()}}\label{fig:tsc_and_adjtimex}
\end{figure}

The main problem is that not only neither $t_0$ nor $r_i$ can be explicitly estimated, but also $r_i$ varies with time as shown in Figure \ref{fig:tsc_drift}.
Thus, current protocols periodically update $s_i^o$ and $x_i^o$ in order to keep track of the changes of $r_i$. These updates are made using the {\it offset} between the current estimate $x_i(t)$ and the global time $t$, i.e. $D_i^x(t)=t-x_i(t)$, and the {\it relative frequency error} that is computed using two offset measurements separated by $\tau$ seconds, i.e.   
\begin{equation}
f^{err}_i(t):=\frac{D_i^x(t)-D_i^x(t-\tau)}{x_i(t)-x_i(t-\tau)} = \frac{1-r_is_i^o}{r_is_i^o}.
\end{equation}
Figure \ref{fig:xt3} provides an illustration of these measurements.



%

\vspace{-.25cm}

To understand the differences between  current protocols, we first rewrite the evolution of $x_i(t)$ based only on the time instants $t_k$ in which the clock corrections are performed. We allow the skew correction $s_i^o$ to vary over time, i.e. $s_i(t_k)$, and write $x_i(t_{k+1})$ as a function of $x_i(t_{k})$. Thus, we obtain
\begin{subequations}
\label{eq:double_int}
\begin{align}
\centering
x_i(t_{k+1}) &= x_i(t_k) +  \tau r_is_i(t_k)  + u_i^x(t_k)\label{eq:double_int_ia}\\
s_i(t_{k+1}) &=s_i(t_k) +  u_i^s(t_k)\label{eq:double_int_ib}
\end{align}
\end{subequations}
where $\tau = t_{k+1}-t_k$ is the time elapsed between adaptations; also known as poll interval~\cite{mills_network_2006}.
The values $u_i^x(t_k)$ and  $u_i^s(t_k)$ represent two different types of corrections that a given protocol chooses to do at time $t_k$ and are usually implemented within the interval $(t_k,t_{k+1})$. $u_i^x(t_k)$ is usually referred to as {\it offset correction } and $u_i^s(t_k)$ as  {\it skew correction}.\footnote{These corrections can be implemented in Linux OS using the {\it adjtimex()} interface to update the system clock or by maintaining a virtual version of $x_i(t)$ and directly applying the corrections to it, as in IBM CCT~\cite{froehlich_achieving_2008} and RADclock~\cite{veitch_robust_2009}. The latter gives more control on how the corrections are implemented since it does not depend on kernel's routines.}
See Figure \ref{fig:adjtimex} for an illustration of their effect on the linux time.
%

 We now proceed to summarize the different types of adaptations implemented by current protocols. The main differences between them are whether they use offset corrections, skew corrections, or both, and whether they update using offset values $D_i^x(t_k)$, relative frequency errors $f_i^{err}(t_{k})$, or both.

\subsection{ Offset corrections}\label{sssec:only_offset}
These corrections consist in keeping the skew fixed and periodically introducing time changes of size  $u_i^x(t_k) = \kp_1 D_i^x(t_k)$ or $u_i^x(t_k) = \kp_1 D_i^x(t_k) + \kp_2 f_i^{err}(t_k)$ 
where $\kp_1,\kp_2>0$.
They are used by NTPv3~\cite{mills_network_1992}  and NTPv4~\cite{mills_network_2006} respectively under ordinary conditions.

\begin{figure}
       \begin{subfigure}[b]{0.49\columnwidth}
               \centering
               \includegraphics[width=.9\columnwidth,height = .5\columnwidth]{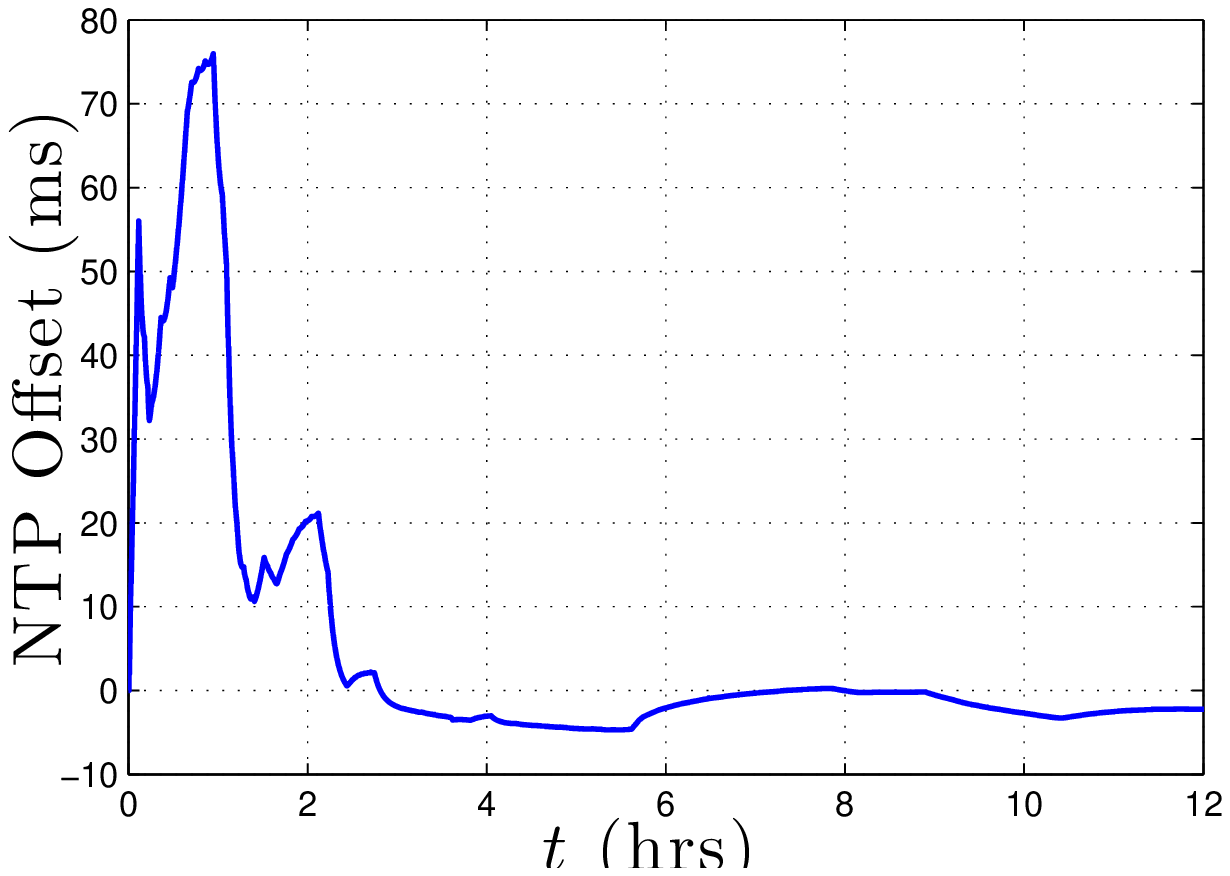}
               \caption{NTP initialization period}\label{fig:ntp_init}
        \end{subfigure}
        \begin{subfigure}[b]{0.49\columnwidth}
               \centering
		\includegraphics[width=.9\columnwidth,height = .5\columnwidth]{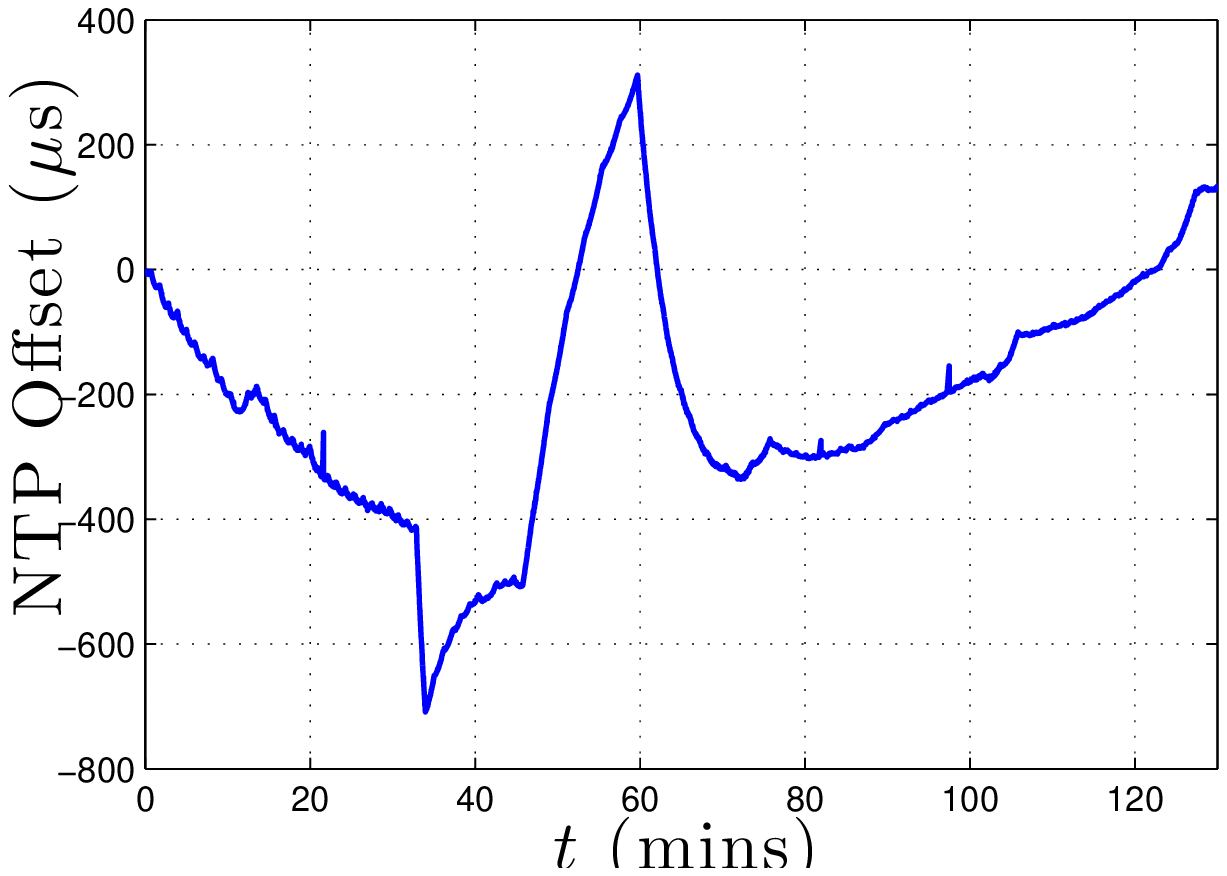}
               \caption{NTP in normal regime}\label{fig:ntp_normal}
        \end{subfigure}%
        \vspace{-.25cm}
        \caption{ Variations of NTP time using TSC as reference }\label{fig:ntp}
\end{figure}

\begin{figure}
       \begin{subfigure}[b]{0.33\columnwidth}
               \centering
               \includegraphics[width=\columnwidth]{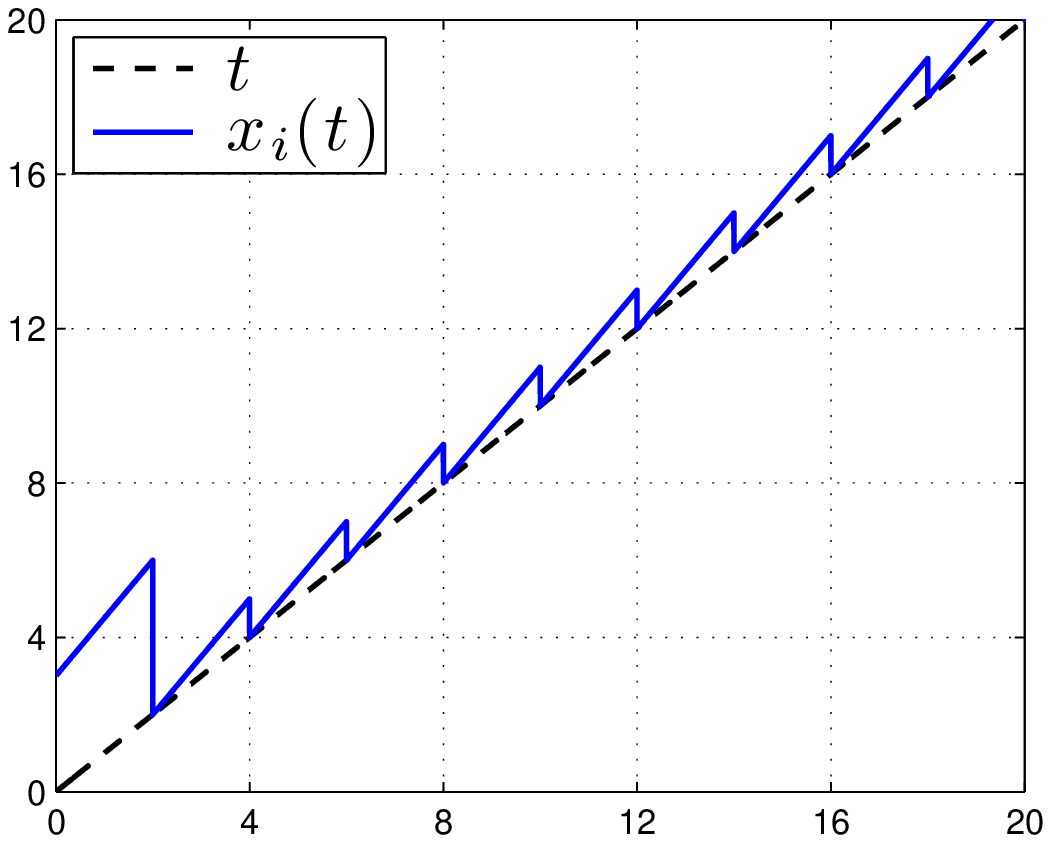}
               \caption{\text{Offset corrections}\quad\text{        }\quad\quad}\label{fig:alg-offset}
        \end{subfigure}
        \begin{subfigure}[b]{0.33\columnwidth}
               	\centering
               	\includegraphics[width=\columnwidth]{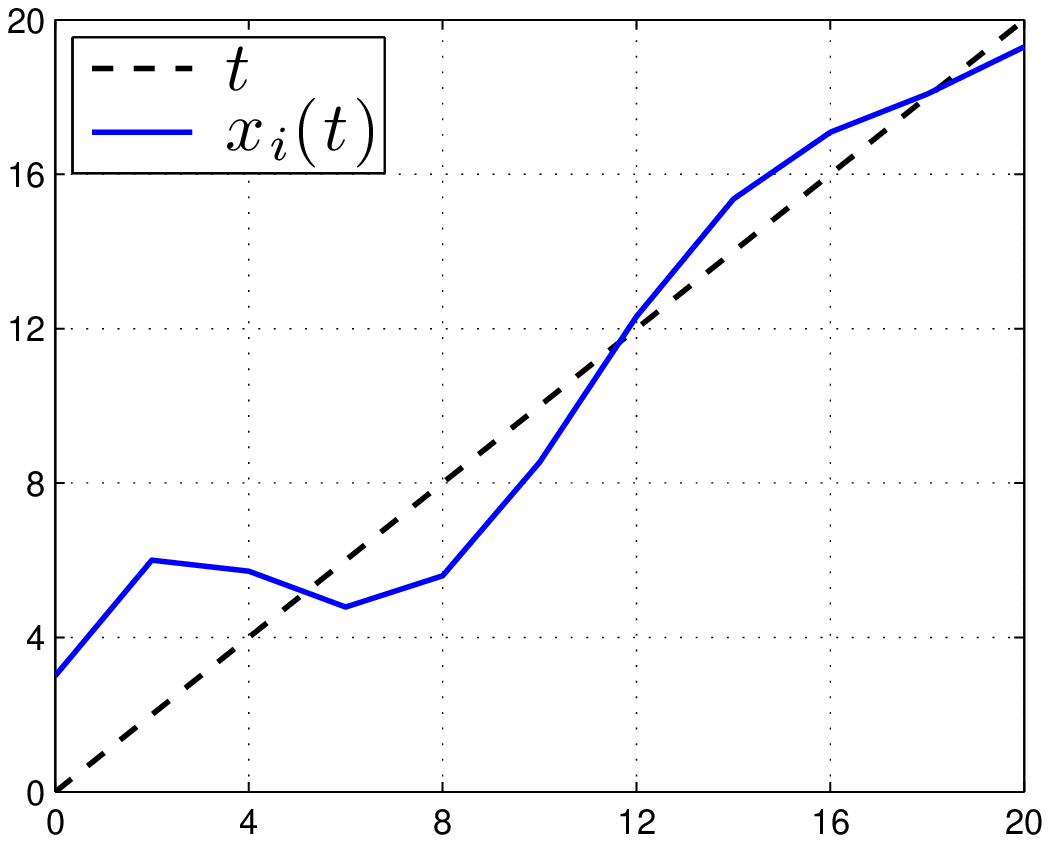}
				\caption{\text{Skew corrections}\quad\text{           }\quad\quad}\label{fig:alg-skew}
        \end{subfigure}%
        \begin{subfigure}[b]{0.33\columnwidth}
               	\centering
               	\includegraphics[width=\columnwidth]{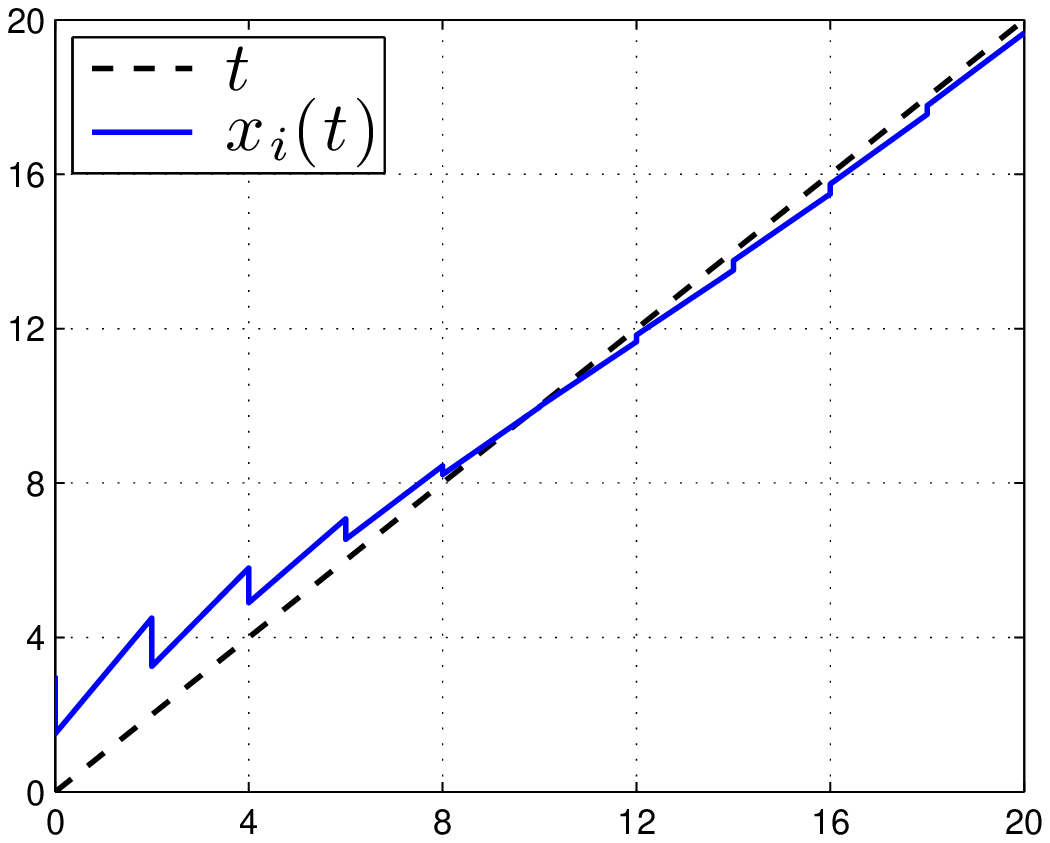}
				\caption{Offset and skew corrections}\label{fig:alg-offsetskew}
        \end{subfigure}%
        \vspace{-.25cm}
        \caption{Current Protocols Adaptation}\label{fig:alg}
\end{figure}

These protocols have in general a slow initialization period as shown in Figure \ref{fig:ntp_init}. This is because the algorithm must first obtain a very  accurate estimate of the initial frequency error $f_i^{err}(t_0)$. 
Furthermore, these updates usually generate non-smooth time evolutions as in Figures \ref{fig:ntp_normal}  and \ref{fig:alg-offset}, and should be done carefully since they might introduce backward jumps ($x_i(t_{k+1})<x_i(t_k)$), which can be problematic for some applications.

\subsection{Skew corrections}\label{sssec:only_skew}
Another alternative that avoids using steep changes in time is proposed by the IBM CCT solution~\cite{froehlich_achieving_2008}. This alternative does not introduce any offset correction, i.e. $u_i^x(t_k) = 0$, and updates the skew $s_i(t_k)$ by 
$u_i^s(t_k) = \kappa_1 D_i^x(t_k) + \kappa_2 f_i^{err}(t_k)$.

The behavior of this algorithm is shown in Figure \ref{fig:alg-skew}. In \cite{xie_consensus_2011} it was shown for a slightly modified version of it 
(used $r_is_i(t_k)f_i^{err}(t_k)$ instead of $f_i^{err}(t_k)$) the algorithm can achieve synchronization for very diverse network architectures.

However, the estimation of $f_i^{err}(t_k)$ is nontrivial as it is constantly changing with subsequent updates of $s_i(t_k)$ and it usually involves sophisticated computations~\cite{zhang_clock_2002,kim_tracking_2012}.



\subsection{Skew and offset corrections}

This type of corrections allow dependence on only offset information $D_i^x(t_k)$ as input to $u_i^x(t_k)$ and $u_i^s(t_k)$. For instance, in \cite{carli_networked_2010} the update
$
u_i^x(t_k) = \kappa_1D_i^x(t_k) \text{ and } 
u_i^s(t_k) = \kappa_2D_i^x(t_k) 
$
was proposed.

This option allows the system to achieve synchronization without any skew estimation. But the cost of achieving it, is introducing  offset corrections in $x_i(t)$ as shown in Figure \ref{fig:alg-offsetskew}. Therefore, it suffers from the same problems discussed in \ref{sssec:only_offset}.

\section{Skewless Network Synchronization}\label{sec:algorithm}

We now present an algorithm that overcomes the limitations of the solutions described in Section \ref{sec:comp_clocks}. In other words, our solution  has the following two properties:
\begin{enumerate}
\item Continuity: The protocol does not introduce steep changes on the time value, i.e. $u_i^x(t_k)\equiv 0$.
\item Skew independence: The protocol does not use skew information $f_i^{err}(t_k)$ as input.
\end{enumerate}
A solution with these properties will therefore prevent unnecessary offset corrections that produce jitter and will be more robust to noise by avoiding skew estimation.
After describing and motivating our algorithm, we show how the updating rule can be  implemented in the context of a network environment.

\begin{figure}
\centering
\includegraphics[width=.85\columnwidth]{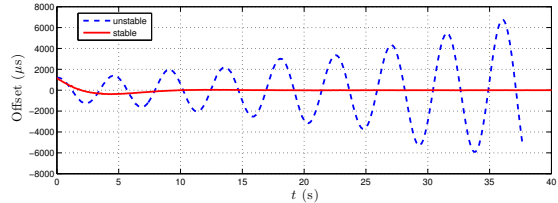}
\caption{Unstable clock steering using only offset information \eqref{eq:double_int_pure} and stable clock steering based on exponential average compensation~\eqref{eq:algorithm} }\label{fig:unstable}
\end{figure}

The motivation behind the proposed solution comes from trying to compensate the problem that arises when one tries to naively impose properties 1) and 2), i.e. using
\begin{align} \label{eq:double_int_pure}
u_i^x(t_k) = 0 \quad \text{ and } \quad u_i^s(t_k) = \kp_1 D_i^x(t_k).
\end{align}
Figure \ref{fig:unstable} shows that this type of clock corrections is unstable; the offset $D_i^x(t_k)$ of the slave clock oscillates with an exponentially increasing amplitude.
\vspace{-.25cm}

The oscillations in Figure \ref{fig:unstable} arise due to the fundamental limitations of using offset to update frequency.
This is better seen in the continuous time version of the system \eqref{eq:double_int} with \eqref{eq:double_int_pure}, i.e.
\[
\dot x_i(t) = r_i s_i(t)\text{ and } \dot s_i(t) = \kp_1 D_i^x(t)
\]
where $\dot x(t) =\frac{d}{dt}x(t)$.
If we consider the offset $D_i^x=t-x_i(t)$ as the system state, then we have 
\[
\dot D_i^x =1-r_is_i\text{ and }
\ddot D_i^x=-\kp_1 r_iD_i^x,
\]
with $\ddot x(t)=\frac{d^2}{dt^2} x(t)$.
 
This is analogous to a spring mass system without friction. Thus, it has two purely imaginary eigenvalues that generate sustained oscillations; see \cite{mallada2008stability,mallada_distributed_2011} for similar examples.\footnote{In the discrete time system the oscillations increase in amplitude since there is a delay between the time the offset is measured $t_k$ and the time the update is made $t_{k+1}$ which makes the system unstable.}
One way to damp these oscillations in the spring-mass case is by adding {\it friction}. This implies  adding a term that includes a frequency mismatch $f_i^{err}(t)$ in our system, which is 
equivalent to the protocols of Section \ref{sssec:only_skew}, and therefore 
undesired. 

However, there are other ways to damp these oscillations using passivity-based techniques from control theory~\cite{ren2008consensus,ren2008distributed}. The basic idea is to introduce an additional state $y_i$ that generates the desired {\it friction} to damp the oscillations.  

Inspired by \cite{ren2008consensus}, we consider the exponentially weighted moving average of the offset
\begin{equation}\label{eq:y}
y_i(t_{k+1}) = pD_i^x(t_k) + (1-p)y_i(t_{k}).
\end{equation}
and update $x_i(t_k)$ and $s_i(t_k)$ using:
%
\begin{align}
u_i^x(t_{k}) = 0 \quad\text{ and }\quad u_i^s(t_{k}) &= \kappa_1 D^x(t_k) -\kappa_2 y(t_{k}).\label{eq:algorithm}
\end{align}
Figure \ref{fig:unstable} shows how the proposed strategy is able to compensate the oscillations without needing to estimate the value of $f_i^{err}(t_k)$.
The stability of the algorithm will depend on how $\kp_1$, $\kp_2$ and $p$ are chosen. A detailed specification of these values is given in Section \ref{ssec:param_vals}.

Finally, since we  are interested in studying the effect of timing loops, we move away from the client-server configuration implicitly assumed in Section \ref{sec:comp_clocks} and allow mutual or cyclic interactions among nodes. 
The interactions between different nodes is described by a  graph $G(V,E)$, where $V$ represents the set of $n$ nodes ($i\in V$) and $E$ the set of {\it directed} edges $ij$; $ij\in E$ means node $i$ can measure its offset with respect to $j$, $D_{ij}^x(t_k) = x_j(t_k)-x_i(t_k)$.

Within this context, a natural extension of  \eqref{eq:y}-\eqref{eq:algorithm} is to substitute $D_i^x(t_k)$ with the weighted average of $i$'s neighbors offsets.
Thus, we propose the following algorithm to update the clocks in the network.

\vspace{.15cm}
\noindent{\bf Algorithm 1 (Alg1):}
For each computer node $i$ in the network, perform the following actions:\noindent
\begin{itemize}[-]
\item  Compute the time offsets ($D^x_{ij}(t_k)$) from $i$ to every neighbor $j$ at time $t_k$.
\item  Update the skew $s_i(t_{k+1})$ and the moving average $y_i(t_{k+1})$ at time $t_{k+1}$ according to:
\begin{subequations}
 \label{eq:system}
 \begin{align}
s_i(t_{k+1}) =&s_i(t_{k}) +  \kp_1\sum_{j\in \mathcal N_i} \alpha_{ij}D^x_{ij}(t_k)
- \kappa_2y_i(t_{k}) \label{eq:system_si}\\
y_i(t_{k+1}) =& p\sum_{j\in \mathcal N_i}\alpha_{ij}D^x_{ij}(t_k) + (1-p)y_i(t_{k})\label{eq:system_yi}
 \end{align}
\end{subequations}
where $\mathcal N_i$ represents the set of neighbors of $i$ and the weights $\alpha_{ij}$ are positive.
\end{itemize}

\vspace{.15cm}
Using this algorithm, many servers can affect the final frequency of the system. Thus, when the system synchronizes, we have
\begin{equation}\label{eq:synchronization}
x_i(t_k) = r^*(t_k-t_0) + x^*\quad i\in V.
\end{equation}
$r^*$ and $x^*$ are possibly different from their ideal values $1$ and $t_0$. Their final values depend on the initial condition of all different clocks as well as the topology, which we assume to be a connected graph in this paper.



\section{Analysis}\label{sec:analysis}
We now analyze the asymptotic behavior of system \eqref{eq:system} and provide a necessary and sufficient condition on the parameter values that guarantee its convergence to \eqref{eq:synchronization}. The techniques used are  drawn from the control literature, e.g. \cite{carli_networked_2010} and \cite{xie_consensus_2011}, yet its application in our case is nontrivial.

\subsubsection*{Notation}
We use $\0_{m\times n}$ ($\1_{m\times n}$) to denote the matrices of all zeros (ones) within $\mathds R^{m\times n}$ and  $\0_{n}$ ($\1_{n}$) to denote the column vectors of appropriate dimensions. $I_n\in\mathds R^{n\times n}$ represents the identity matrix.  Given a matrix $A\in\mathds R^{n\times n}$ with Jordan normal form $A=PJP^{-1}$, let $n_A\leq n$ denote the total number of Jordan blocks $J_l$ with $l\in\mathcal I(A):=\{1,...,n_A\}$. 
We use $\mu_l(A)$, $l\in\{1,\dots,n\}$ or just $\mu(A)$ to denote the eigenvalues of $A$, and order them decreasingly $|\mu_{1}(A)|\geq\dots\geq |\mu_{n}(A)|$. Finally, $A^T$ is the transpose of $A$, $A_{ij}$ is the element of the $i$th row and $j$th column of $A$ and $a_i$ is the $i$th element of the column vector $a$, i.e. $a=[a_i]^T$. 

It is more convenient for the analysis to use a vector form representation of \eqref{eq:system} given by
\begin{equation}\label{eq:system_z}
z_{k+1} = A z_k
\end{equation}
where $z_k:=[ x(t_{k})^T s(t_{k})^T y(t_{k})^T]^T\in\mathds R^{3n}$, 
\begin{equation*}
A:=
\left[
\begin{array}{ccc}
I_n & \tau R &\0_{n\times n} \\
-\kappa_1L & I_n & -\kappa_2I_n \\
p(-L) & \0_{n\times n} & (1-p) I_n
\end{array}
\right]\in \mathds{R}^{3n\times 3n},
\end{equation*}
$R\in\mathds R^{n\times n}$ is the diagonal matrix with elements $r_i$ and $L\in\mathds R^{n\times n}$ is the Laplacian matrix associated with $G(V,E)$,
\begin{align*}
L_{ii} =\alpha_{ii}:= \sum_{j\in\mathcal N_i} \alpha_{ij} \text{ and } L_{ij} &=
\begin{dcases*}
 -\alpha_{ij}                        & if $ij \in E$,\\
0                      & otherwise.
\end{dcases*}
\end{align*}

The convergence analysis of this section is done in two stages. First, we provide necessary and sufficient conditions for synchronization in terms of the eigenvalues of $A$ (Section \ref{ssec:asymptotic_behavior}) and then use Hermite-Biehler Theorem \cite{bhattacharyya_robust_1995} to relate these eigenvalues with the parameter values  that can be directly used in practice (Section \ref{ssec:param_vals}). All the proof details are included in the appendix for interested readers.

\subsection{Asymptotic Behavior}\label{ssec:asymptotic_behavior}
We start  by studying the asymptotic behavior of \eqref{eq:system_z}. That is, we are interested in finding under what conditions the series of elements $\{x_i(t_k)\}$ converge to \eqref{eq:synchronization} as $t_k$ goes to infinity.


Consider the Jordan normal form \cite{horn_matrix_1990} of 
\[
A := \left[ \zeta_1 \quad ... \quad\zeta_{3n}\right] 
J
\left[ \eta_1 \quad ... \quad\eta_{3n}\right]^T
\]
 where $J=\blockdiag(J_l)_{l\in \mathcal I(A)}$, $\zeta_i$ and $\eta_i$ are the right and left generalized eigenvectors of $A$ such that
\begin{align*}
\zeta_i^T\eta_j =
\begin{dcases*}
1 & if $j=i$,\\
0 & otherwise.
\end{dcases*}
\end{align*}
The crux of the analysis comes from understanding the relationship between the multiplicity of the eigenvalue $\mu(A)=1$ and the eigenvalue $\mu(L) = 0$,  and their corresponding eigenvectors. This is captured in the next two lemmas.

\begin{lemma}[Eigenvalues of $A$ and Multiplicity of $\mu(A) = 1$] \label{lem:multiplicity}
$A$ has an eigenvalue $\mu(A)=1$ with multiplicity $2$ if and only if the graph $G(V,E)$ is connected, $\kappa_1\neq\kappa_2$ and $p>0$.

Furthermore, $\mu_l(A)$ are the roots of
\begin{equation}\label{eq:g_l}
g_l(\lambda):=(\lambda-1)^2(\lambda-1+p)+  [(\lambda-1)\kappa_1+\kappa_2-\kappa_1]\nu_l
\end{equation}
where $\nu_l = \mu_l(\tau LR)$ and satisfies
\begin{equation}\label{eq:nu_condition}
\nu_n=0<|\nu_l|\text{ for }l\in\{1,\dots,n-1\}.
\end{equation}
\end{lemma}


\begin{lemma}[Jordan Chains of $\mu(A)=1$ and $\mu(A)=1-p$]\label{lem:jordan_chains}
Under the conditions of Lemma \ref{lem:multiplicity} the right and left Jordan chains, $(\zeta_1,\zeta_2)$ and $(\eta_2,\eta_1)$ respectively, associated with $\mu(A)=1$ and the eigenvectors $\zeta_3$ and $\eta_3$ associated with $\mu(A)=1-p$ are given by 
\begin{equation}\label{eq:zeta}
[\zeta_1\;\zeta_2\;\zeta_3]= \left[\begin{array}{ccc}
\1_n & \1_n                                & -\frac{\tau\kp_2}{p^2}\1_n \\ 
\0_n & \frac{(R^{-1}\1_n)}{\tau} & \frac{\kp_2}{p}R^{-1}\1_n \\ 
\0_n & \0_n                                & R^{-1}\1_n 
\end{array}
\right] \text{ and }
\end{equation}
\begin{equation}\label{eq:eta}
[\eta_1\;\eta_2\;\eta_3]=\gamma \left[\begin{array}{ccc}
 R^{-1}\xi 																&	\0_n							& 	\0_n		\\ 
 -\tau\xi 																&	\xi								&	\0_n		\\ 
 \tau\kp_2(\frac{1}{p}+\frac{1}{p^2})\xi 	& -\frac{\kp_2}{p}\xi 		&	\xi
\end{array}
\right]
\end{equation}
where $\xi$ is the unique normalized left eigenvector of $\mu(L)=0$ ($\sum_{i=1}^n \xi_i=1$) and  $\gamma$ is the $\xi_i$-weighted harmonic mean of $r_i$, i.e. $
\frac{1}{\gamma} = \1_n^T R^{-1}\xi = \sum_{i=1}^{n} \frac{\xi_i}{r_i}.
$
\end{lemma}


The proof of Lemmas \ref{lem:multiplicity} and \ref{lem:jordan_chains} can be found in the Appendices \ref{app:lemma1} and \ref{app:lemma2}.
We now proceed to state our main convergence result.

\begin{theorem}[Convergence]\label{th:convergence}
The algorithm \eqref{eq:system_z} achieves synchronization for any initial conditions if and only if the graph $G(V,E)$ is connected, $\kappa_1\neq\kappa_2$, $p>0$ and $|\mu_l(A)|<1$ whenever $\mu_l(A)\neq 1$.
Moreover, whenever the system synchronizes, we have
\begin{subequations}\label{eq:xs&rs}
\begin{eqnarray}\label{eq:xs}
x^* = \gamma\sum_{i=1}^n \xi_i\left(\frac{1}{r_i} x_{i}(t_0)+\tau\frac{\kp_2}{p^2}y_i(t_0)\right), \text { and }
\end{eqnarray}
\vspace{-.5cm}
\begin{eqnarray}\label{eq:rs}
r^* = \gamma\sum_{i=1}^{n} \xi_i(s_{i}(t_0) - \frac{\kappa_2}{p}y_{i}(t_0)).
\end{eqnarray}
\end{subequations}
\end{theorem}

Theorem \ref{th:convergence} provides an analytical tool to understand the influence of the different nodes of the graph in the final offset $x^*$ and frequency $r^*$. For example, suppose that we know that node 1 has perfect knowledge of its own frequency ($r_1$) and the UTC time at $t=t_0$ ($x_1(t_0)=t_0$), and configure the network such that node 1 is the {\bf unique leader} like the top node in Figures \ref{fig:topologies}a and \ref{fig:topologies}c. It is easy to show that $\xi_1=1$ and $\xi_i=0$ $\forall i\neq1$.  Then, using \eqref{eq:xs}-\eqref{eq:rs} and definition of $\gamma$ we can see that $\gamma=r_1$ and
\[
x^* = x_1(t_0) + r_1\tau \frac{\kp_2}{p^2}y_1(t_0)\text{ and }
r^* = r_1s_{1}(t_0) - \frac{r_1\kappa_2}{p}y_{1}(t_0).
\]
However, since node 1 knows $r_1$ and $t_0$, it can choose $x_1(t_0)=t_0$, $s_1(t_0)=\frac{1}{r_1}$ and $y_1(t_0)=0$. Thus, we obtain $x^*=t_0$ and $r^*=1$ which implies by \eqref{eq:synchronization} that every node in the network will end up with $x_i(t) = t $. In other words, Theorem \ref{th:convergence} allows us to understand how the information propagates and how we can guarantee that every server will converge to the desired time. 
Notice that the initial condition used for server 1 is equivalent to assuming that server 1 is a reliable source of UTC like an atomic clock for instance.

\subsection{Necessary and sufficient conditions for synchronization}\label{ssec:param_vals}

\begin{figure}
\centering
\includegraphics[width=\columnwidth]{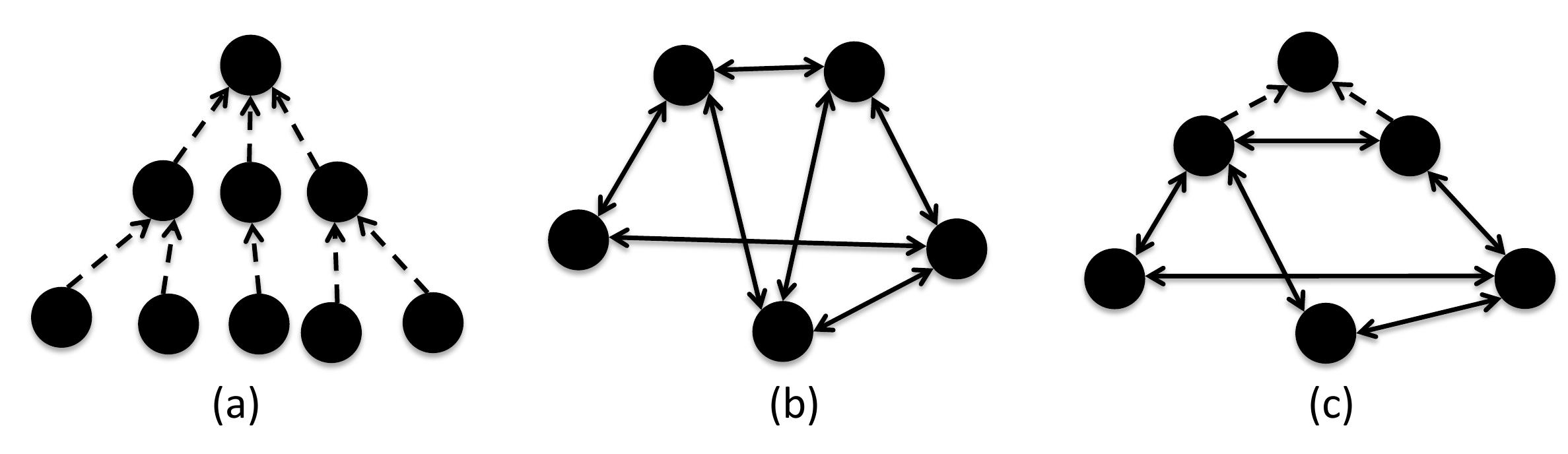}
\caption{Graphs with real eigenvalue Laplacians }\label{fig:topologies}
\end{figure}

We now provide necessary and sufficient conditions in terms of explicit parameter values ($\kp_1$, $\kp_2$ ,$\tau$ and $p$) for Theorem~\ref{th:convergence} to hold.
We will restrict our attention to graphs that have Laplacian matrices with real eigenvalues.
This includes for example trees (Figure \ref{fig:topologies}a), symmetric graphs with $\alpha_{ij}=\alpha_{ji}$ (Figure \ref{fig:topologies}b) and symmetric graphs with a leader (Figure \ref{fig:topologies}c).


\vspace{-.25cm}

The proof consists on studying the Schur stability of $g_l(\lambda)$ and has several steps. We first perform a change of variable that maps the unit circle onto the left half-plane. This transforms the problem of studying the Schur stability into a Hurwitz stability problem which is solved using  Hermite-Biehler Theorem.
\begin{theorem}[Hurwitz Stability (Hermite-Biehler)~\cite{bhattacharyya_robust_1995}]\label{th:HB}
Given the polynomial $P(s)=a_ns^n +...+a_0$, let $P^r(\omega)$ and $P^i(\omega)$ be the real and imaginary part of $P(j\omega)$, i.e. $P(j\omega) = P^r(\omega) +jP^i(\omega)$. Then $P(s)$ is a Hurwitz polynomial if and only if
\begin{enumerate}
\item $a_na_{n-1}>0$ and 2) 
\item The zeros of $P^r(\omega)$ and $P^i(\omega)$ are all simple and real and interlace as $\omega$ runs from $-\infty$ to $+\infty$.
\end{enumerate}
\end{theorem}

We now determine the proper parameter values that guarantee synchronization.
\begin{theorem}[Parameter Values for Synchronization]\label{th:param_sync}
Given a connected graph $G(V,E)$ such that the corresponding Laplacian matrix $L$ has real eigenvalues. The system \eqref{eq:system_z} achieves synchronization 
if and only if
\begin{enumerate}[(i)]
\item $|1-p| <1$ or equivalently $2>p>0$
\item $\frac{2\kp_1}{3p} >\kappa_1-\kappa_2 >0$ and (iii) 
$\tau< \frac{p(\kappa_2 -p(\kappa_1-\kappa_2))}{\mu_{\max}(\kappa_1-p(\kappa_1-\kappa_2))^2}$
\end{enumerate}
where $\mu_{\max}$ is the largest eigenvalue of $LR$.
\end{theorem}

Even though $\mu_{\max}$ depends on $r_i$ which is in general unknown, it is easy to show that $\mu_l(LR)\leq\hat r_{\max}\mu_l(L)$  where $\hat r_{\max}$ is an upper bound of the maximum rate deviation $r_i$. Furthermore, using Greshgorin's circle theorem, it is easy to show that $\mu_{\max}(L)\leq2\alpha_{\max} :=2 \max_i \alpha_{ii}$.
Therefore, if we set
\begin{equation}\label{eq:tau_bound}
\tau <  \frac{p(\kp_2 - \delta\kp p)}{ 2\alpha_{\max}\hat r_{\max}(\kp_1-\delta\kp p)^2}
\end{equation}
convergence is guaranteed for {\bf every connected graph} with real  eigenvalues.

\section{Experiments}\label{sec:experimental_results} 

To test our solution and analysis, we implement  an asynchronous version of Algorithm 1 (Alg1) in C using the IBM CCT solution as our code base.
Every node perform its own measurements and updates every $\tau$ seconds using \eqref{eq:system}, but not necessarily at the same instants $t_k$.

Our program reads the TSC counter directly using the {\verb rdtsc } assembly instruction to minimize reading latencies and maintains a virtual clock that can be directly updated. The list of neighbors is read from a configuration file and whenever there is no neighbor, the program follows the local Linux clock.  Finally, offset measurements are taken using an improved ping pong mechanism proposed in~\cite{froehlich_achieving_2008}. 

\begin{figure}
\centering
\includegraphics[width=.8\columnwidth]{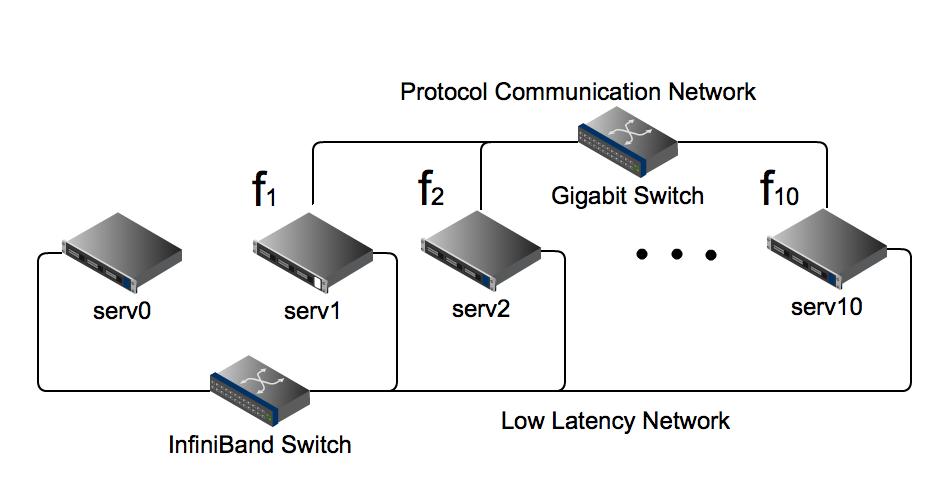}
\caption{Testbed of IBM BladeCenter blade servers} \label{fig:testbed}
\end{figure}

We run our skewless protocol in a cluster of IBM BladeCenter LS21 servers with two AMD Opteron processors of $2.40$GHz, and 16GB of memory. As shown in Figure \ref{fig:testbed}, the servers serv1-serv10 are used to run the protocol. The offset measurements are taken through a Gigabit Ethernet switch. Server serv0 is used as a reference node and gathers time information from the different nodes using a Cisco 4x InfiniBand Switch that supports up to 10Gbps between any two ports and up to 240Gbps of aggregate bandwidth. This minimizes the error induced by the data collecting process.

We use this testbed to validate the analysis in Section \ref{sec:analysis}. 
First, we illustrate the effect of different parameters and analyze
the effect of the network configuration on convergence (Experiment 1).
Then we present a series of configurations that demonstrate
how connectivity between clients is useful in reducing the
jitter of a noisy clock source (Experiment 2). And finally, we compare the
performance of our protocol with respect to NTP version 4 (Experiment 3) and 
IBM CCT (Experiment 4).

We will use several performance metrics to evaluate Alg1. For instance, the {\it mean relative deviation } from the leader which is defined as the root mean square of the node's offset with respect to the leader, i.e.  $\sqrt{S_n}$ with
\begin{equation}
S_n = \frac{1}{{n-1}}{\sum_{i=2}^{n}\left<(x_i-x_1)^2\right>},
\end{equation}
where $<\cdot>$ amounts to the sample average.
We will also use the $99\%$ Confidence Interval $CI_{99}$ and the maximum offset ($CI_{100}$) as metrics of accuracy. For example, if $CI_{99}=10\mu s$, then  $99\%$ of the offset samples will be within 10$\mu s$ of the leader. 

Unless explicitly stated, the default parameter values are
\begin{align}\label{eq:values1}
p=0.99,\quad \kp_1 = 1.1, \quad  \kp_2 = 1.0 \text{ and }
\alpha_{ij}=\frac{c}{|\mathcal N_i|}.
\end{align}
The scalar $c$ is a commit or gain factor that will allow us to compensate the effect of $\tau$. Notice that by definition of $\alpha_{ij}$, $\alpha_{ii}=c$ for every node that is not the leader.


Moreover, these values immediately satisfy (i) and (ii) of Theorem \ref{th:param_sync} since $1-p=0.01$ and $\frac{2\kp_1}{3p} = 0.7407 > \kp_1-\kp_2 = 0.1$. The remaining condition can be satisfied by modifying $\tau$ or equivalently $c$. Here, we choose to fix $c=0.7$ which makes condition (iii)
\[
\tau< \frac{890.1}{\mu_{\max}}\text{ms}.
\]
For fixed polling interval $\tau$, the stability of the system depends on the value of $\mu_{\max}$, which is determined by the underlying network topology and the values of $\alpha_{ij}$.

\begin{figure}[htp]
\centering
\includegraphics[width=.55\columnwidth]{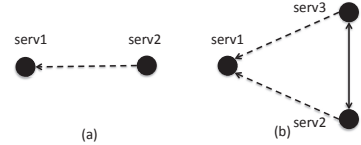}
\caption{Effect of topology on convergence: (a) Client-server configuration;  (b) Two clients connected to server and mutually connected.}\label{fig:23nodes}
\vspace{-.25cm}
\end{figure}

\begin{figure}[htp]
        \begin{subfigure}[b]{0.49\columnwidth}
               \centering
               \includegraphics[width=\columnwidth,]{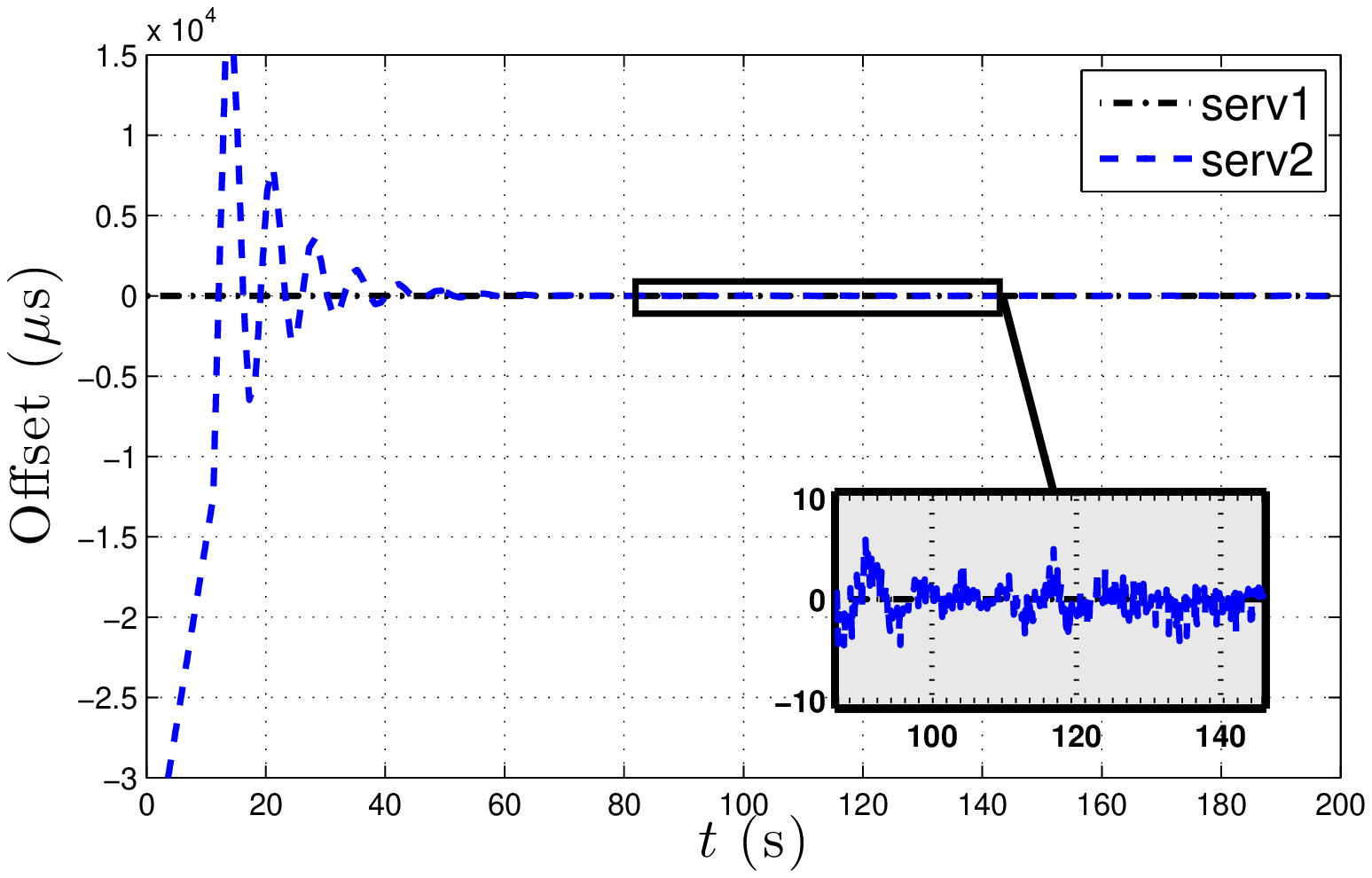}
               \caption{Client server configuration with $\tau = 1$s. The client converges and the algorithm is stable.}\label{fig:1000ms_21}
        \end{subfigure}%
        ~ 
        \begin{subfigure}[b]{0.49\columnwidth}
               \centering
               \includegraphics[width=\columnwidth]{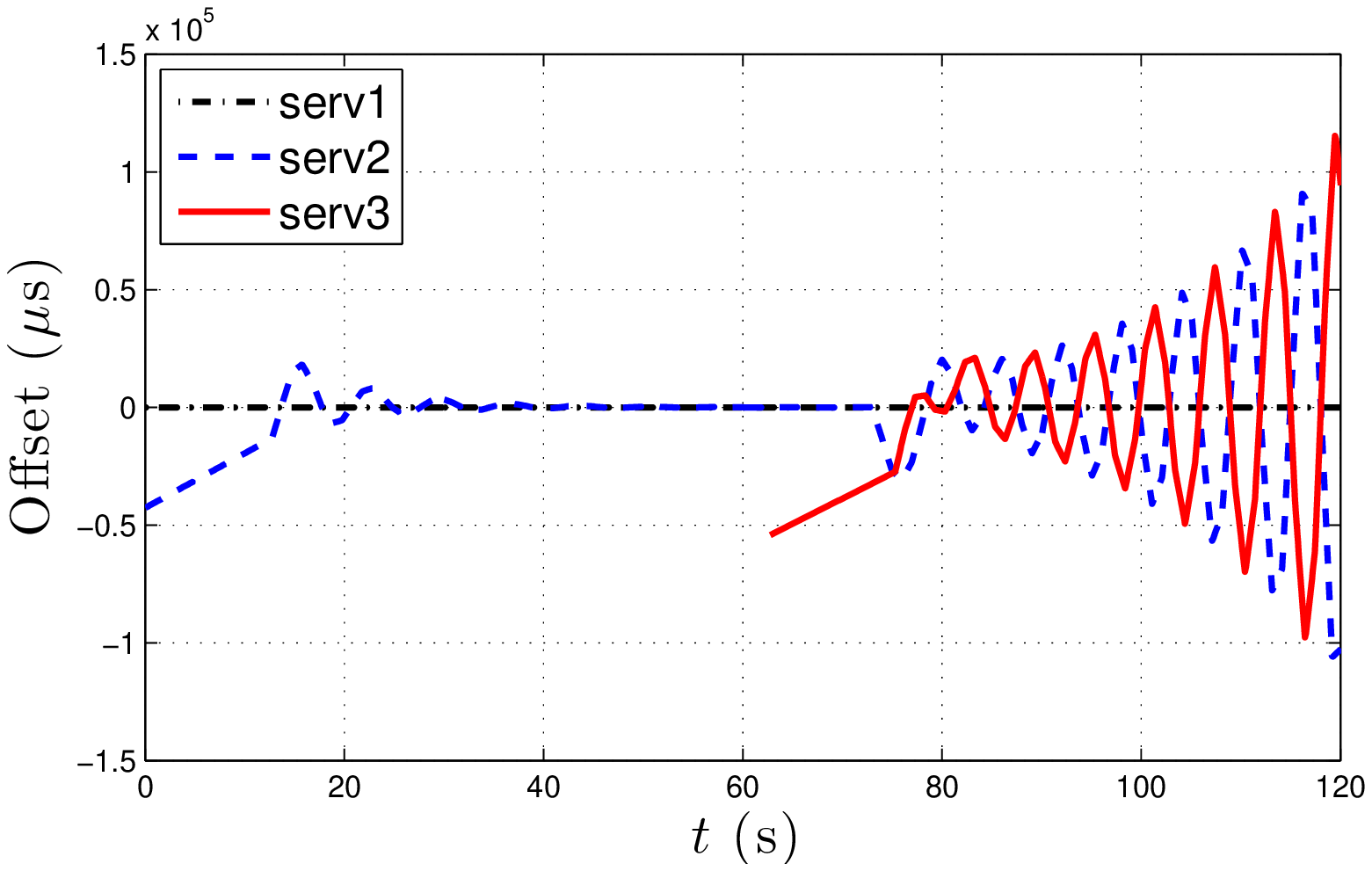}
               \caption{Two clients mutually connected with $\tau = 1$s. The algorithm becomes unstable.} \label{fig:1000ms_213}
        \end{subfigure}
        \caption{Loss of stability by change in the network topology}\label{fig:offsets}
\end{figure}

\vspace{.15cm}
\noindent{\bf Experiment 1 (Convergence):}
We first consider the client server configuration described in Figure \ref{fig:23nodes}a with a time step
$
\tau = 1\text{s}.
$ 
 In this configuration $\mu_{\max} \approx c=0.7$ and therefore condition (iii) becomes $\tau<1.2717$s. Figure \ref{fig:1000ms_21} shows the offset between serv1 (the leader) and serv2 (the client) in microseconds.
 There we can see how serv2 gradually updates $s_2(t_k)$ until the offset becomes insignificant.

Figure \ref{fig:1000ms_21} tends to suggest that the set of parameters given by \eqref{eq:values1} and $\tau=1$s are suitable for deployment on the servers. This is in fact true provided that network is a directed tree as in Figure \ref{fig:topologies}a. 
The  intuition behind this fact is that in a tree, each client connects only to one server. Thus, those connected to the leader will synchronize first and then subsequent layers will follow.

However, once loops appear in the network, there is no longer a clear dependency since two given nodes can mutually get information from each other. This type of dependency might make the algorithm unstable.
Figure \ref{fig:1000ms_213} shows an experiment with the same configuration as Figure \ref{fig:1000ms_21} in which serv2 synchronizes with serv1 until a third server (serv3) appears after $60$s. At that moment the system is reconfigured to have the topology of Figure \ref{fig:23nodes}b introducing a timing loop between serv2 and serv3. This timing loop makes the system unstable.

The instability arises since after serv3 starts, the new topology has  $\mu_{\max}\approx 1.5c=1.05$. Thus, the time step condition (iii) becomes
$
\tau< 847.8\text{ms}
$
which is no longer satisfied by $\tau=1$s. 

\begin{figure}[htp]
\centering
\includegraphics[width=.85\columnwidth]{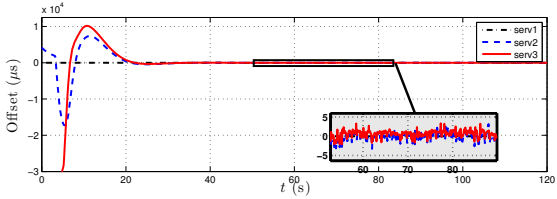}
\caption{Two clients mutually connected with $\tau = 500$ms}\label{fig:500ms_213}
\end{figure}

This may be solved for the new topology (Figure \ref{fig:23nodes}b) by using any $\tau$ smaller than $847.8$ms. However, if we want a set of parameters that is independent of the topology, we can use  \eqref{eq:tau_bound} and notice that $\alpha_{\max}=c$ and $\hat r_{\max}\approx 1$. We choose
\[
\tau = 500\text{ms} < \frac{890.2}{2\alpha_{\max}}\text{ms}=\frac{890.2}{2c}\text{ms}= 635.9\text{ms}.
\]
Figure \ref{fig:500ms_213} shows how  now serv2 and serv3 can synchronize with serv1 after introducing this change.
\vspace{-.25cm}

\vspace{-.25cm}
\noindent
{\bf Experiment 2 (Timing Loops Effect):}
 We now show how timing loops can be used to collectively outperform individual clients when the time source is noisy.

\begin{figure}[htp]
\centering
\includegraphics[width=.8\columnwidth]{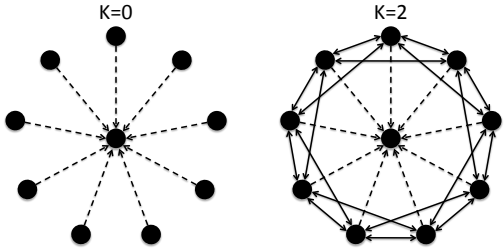}
\caption{Leader topologies with $2K$ neighbors connection. Connections to the leader (serv1) are unidirectional while the connections among clients (serv2 trhough serv10) are bidirectional}\label{fig:wheel}
\end{figure}

We run Alg1 on 10 servers (serv1 through serv10). The connection setup is described in Figure \ref{fig:wheel}. Every node is directly connected unidirectionally to the leader (serv1) and bidirectionally to $2K$ additional neighbors. 
\vspace{-.25cm}

When $K=0$ then the network reduces to a star topology and when $K=4$ the servers serv2 through serv10 form a complete graph.

The dashed arrows in Figure \ref{fig:wheel} show the connections where jitter was introduced. To emulate a link with jitter we added random noise $\eta$ with values taken uniformly from $\{0,1,...,\text{Jitter}_{\max}\}$ on both direction of the communication,
\begin{equation}
\eta \in \{0,1,...,\text{Jitter}_{\max}\}\text{ms}.
\end{equation}

Notice that the arrow only shows a dependency relationship, the ping pong mechanism sends packets in both direction of the physical communication. We used a value of Jitter$_{\max}=10$ms.
Since the error was introduced in both directions of the ping pong, this is equivalent to a standard deviation of $6.05$ms.

 \begin{figure}[htp]
        \begin{subfigure}[b]{0.49\columnwidth}
                \centering
                \includegraphics[width=\columnwidth]{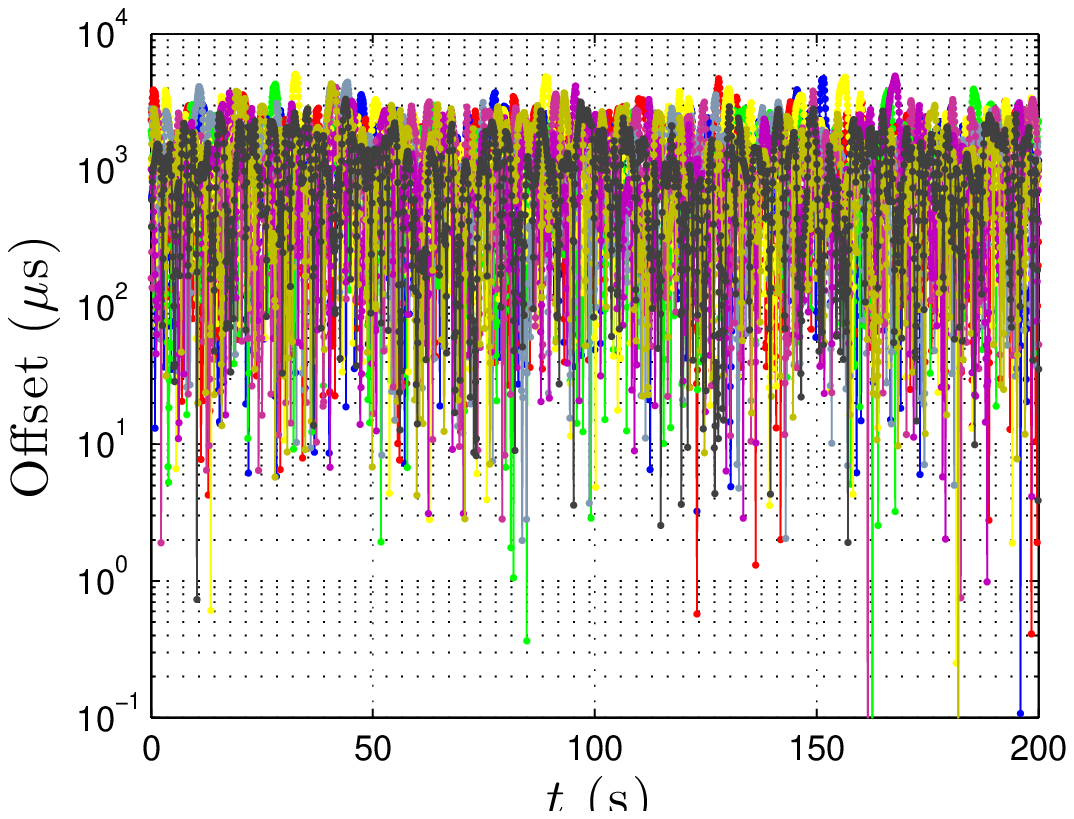}
                \caption{Star topology ($K=0$)}
                \label{fig:offset_star}
        \end{subfigure}%
        ~ 
        \begin{subfigure}[b]{0.49\columnwidth}
                \centering
                \includegraphics[width=\columnwidth]{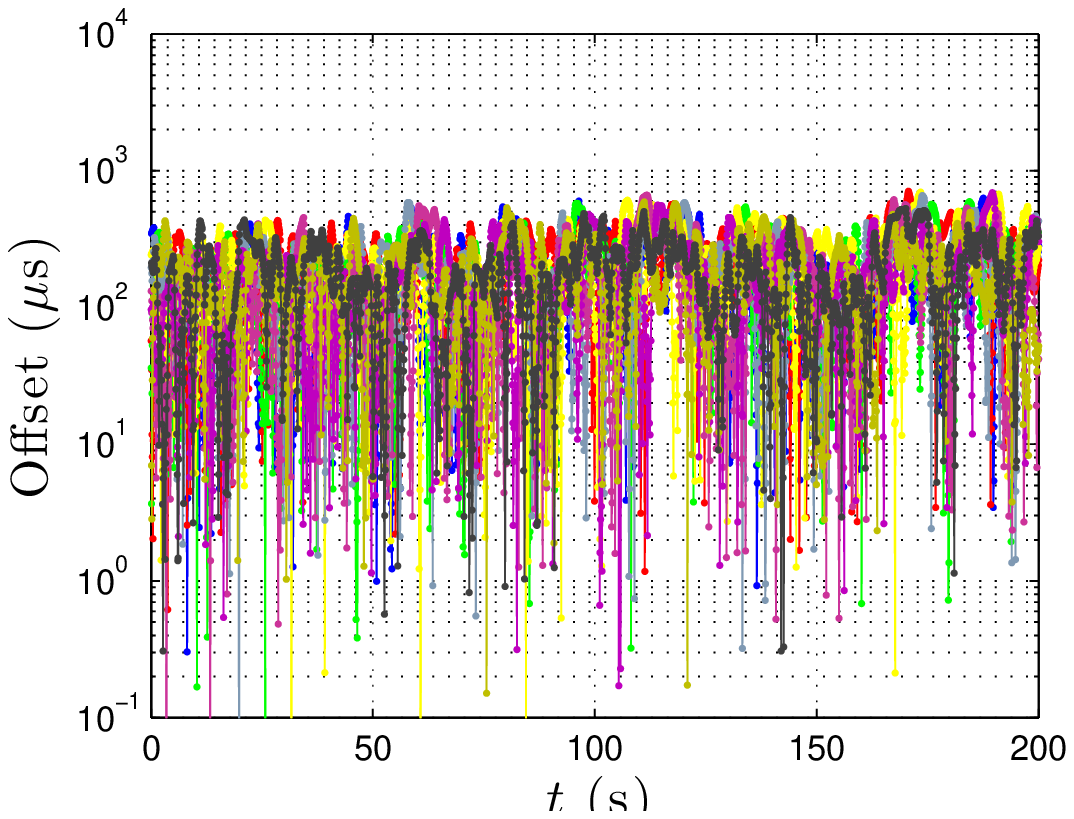}
                \caption{Complete subgraph ($K=4$)}
                \label{fig:offset_K}
        \end{subfigure}
        \vspace{-.25cm}
        \caption{Offset of the nine servers connected to a noisy clock source}\label{fig:offsets}
\end{figure}\vspace{-.25cm}

Figure \ref{fig:offsets} illustrates the relative offset between the two extreme cases; The star topology ($K=0$) is shown in Figure \ref{fig:offset_star}, and the complete subgraph ($K=4$) is shown in Figure \ref{fig:offset_K}.

The worst case offset for $K=0$ is $CI_{100}=5.1$ms which is on the order of the standard deviation of the jitter. However, when $K=4$ we obtain a worst case offset of $CI_{100}=690.8\mu$s, an order of magnitude improvement. 

\begin{figure}[htp]
\centering
\includegraphics[width=.85\columnwidth]{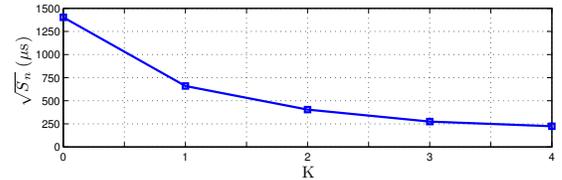}
\caption{Effect of the client's communication topology on the mean relative deviation. As the connectivity increases ($K$ increases) the mean relative deviation is reduced by factor of  $6.26$, i.e. a noise reduction of approx. 8dB.}\label{fig:error_wheel}
\end{figure}\vspace{-.25cm}

The change on the mean relative deviation $\sqrt{S_n}$ as the connectivity among clients increases from isolated nodes ($K=0$) to a complete subgraph ($K=4$) is studied in Figure \ref{fig:error_wheel}.
The results presented show that even without any offset filtering mechanism the network itself is able to perform a distributed filtering that achieves an improvement of up to a factor of $6.26$ or equivalently a noise reduction of almost 8dB.

\vspace{.15cm}
\noindent{\bf Experiment 3 (Comparison with NTPv4):}
We now perform a thorough comparison between our protocol (Alg1) and NTPv4.
We will use the one hop configuration of Figure \ref{fig:23nodes}b but without the bidirectional link. Here,  server serv1 is set as NTP server and as leader of Alg1, server serv2 has a client running NTPv4 and server serv3 a client running our protocol.

In order to make a fair comparison, we need both algorithms to use the same polling interval. Thus, we fix $\tau=16$sec. This can be done for NTP by setting the parameters  {\verb minpoll } and {\verb maxpoll } to $4$ ($2^4=16$secs). 
The remainder parameter values for Alg1 are given by
\begin{align}\label{eq:values3}
p=1.98,\quad \kp_1 = 1.388 \text{ and }  \kp_2 = 1.374.
\end{align}

\begin{figure}[htp]
       \begin{subfigure}[b]{0.49\columnwidth}
               \centering
               \includegraphics[width=\columnwidth]{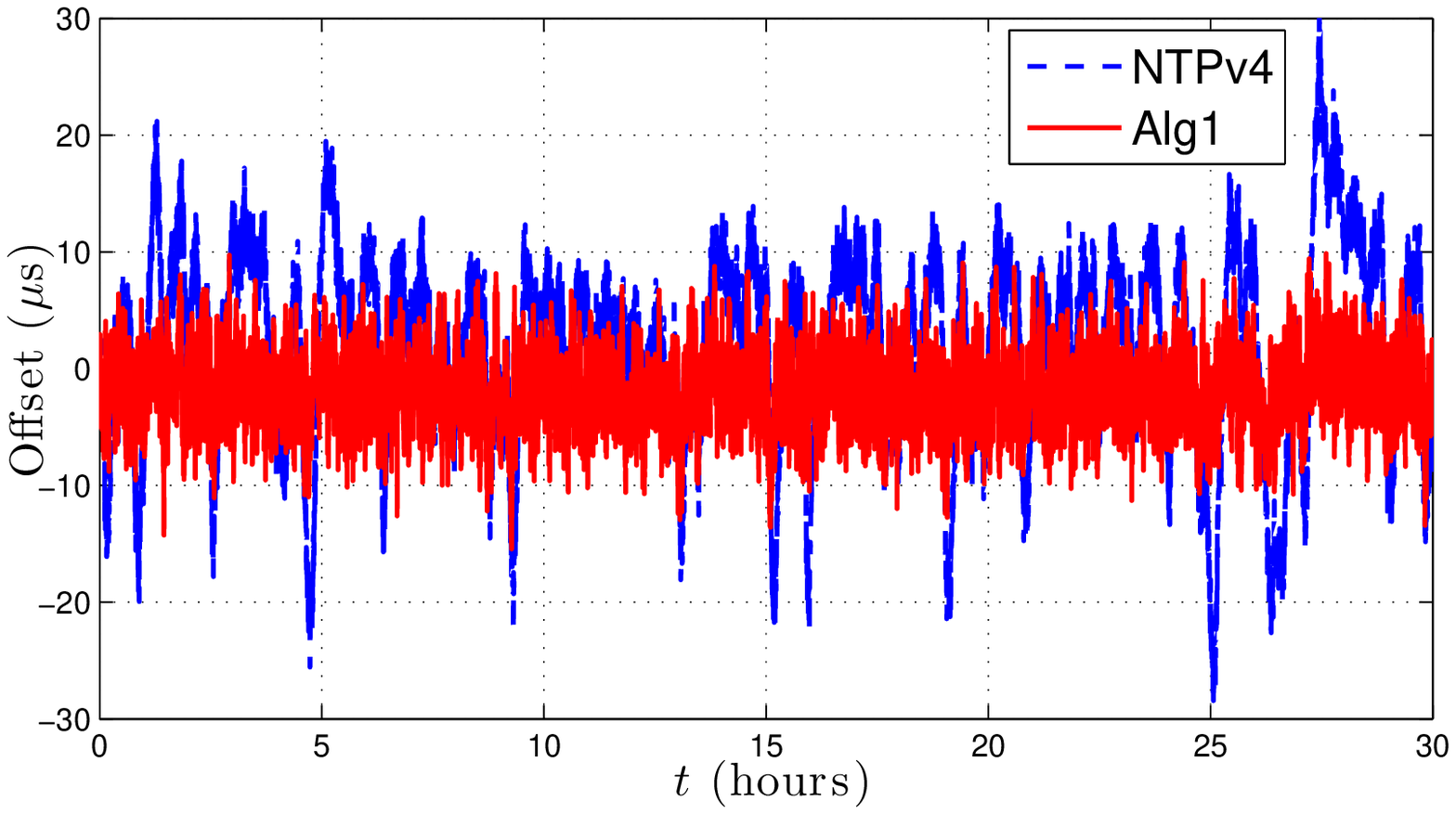}
               \caption{Offset values of NTPv4 and Alg1 for a period of 30 hours.}\label{fig:offset_ntpvsalg1}
        \end{subfigure}
        \begin{subfigure}[b]{0.49\columnwidth}
               \centering
               \includegraphics[width=\columnwidth]{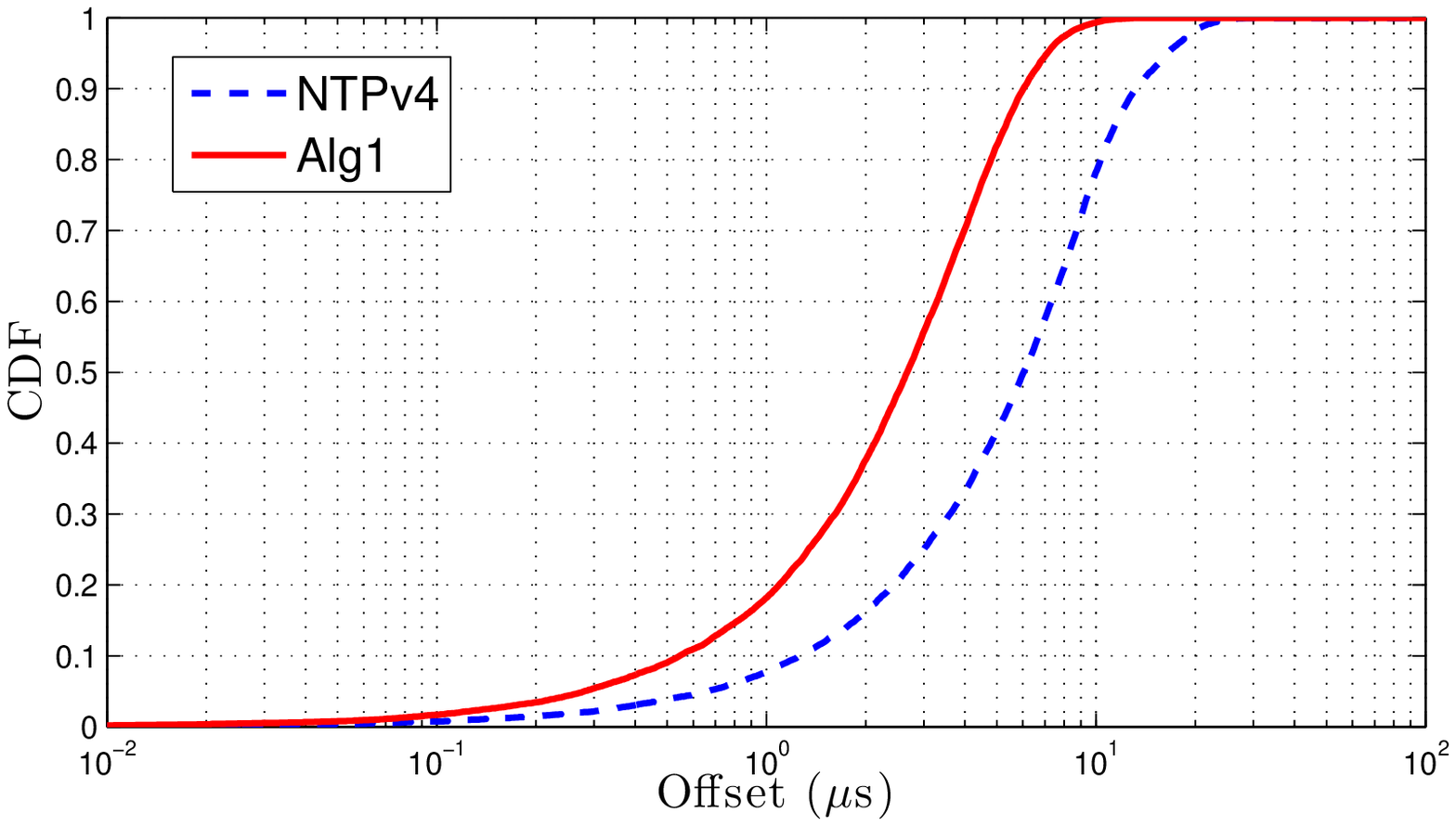}
               \caption{Cummulative Distribution Function}\label{fig:cdf_ntpvsalg1}
        \end{subfigure}%
        \vspace{-.25cm}
        \caption{Performance evaluation between our solution (Alg1) and NTPv4}\label{fig:comparison_ntp}
\end{figure}

Figure \ref{fig:offset_ntpvsalg1} shows the time differences between the clients running NTPv4 and Alg1 (serv2 and serv3) , and the leader (serv1) over a period of 30 hours.  It can be seen that Alg1 is able to track serv1's clock keeping an offset smaller than 10$\mu$s for most of the time while NTPv4 incurs in larger offsets during the same period of time. This difference is produced by the fact that Alg1 is able to react more rapidly to frequency changes while NTPv4 incurs in more offset corrections that generate larger jitter.

\vspace{-.25cm}


A more detailed and comprehensive analysis is presented in Figure \ref{fig:cdf_ntpvsalg1} where we plot the Cumulative Distribution Function (CDF) of the offset samples. That is, the fraction of samples whose time offset is smaller than a specific value. Using Figure \ref{fig:cdf_ntpvsalg1} we compute the corresponding $99\%$ confidence intervals ($CI_{99}$)


Alg1 achieves a performance of $\sqrt{S_n}=3.1\mu$s, $CI_{99}=9.5\mu$s and a maximum offset of $CI_{100}=15.9\mu$s, while NTPv4 obtains
$\sqrt{S_n}=8.1\mu$s, $CI_{99}=21.8\mu$s and a maximum offset of $CI_{100}=28.0\mu$s.
Thus, not only Alg1 achieves a reduction of $\sqrt{S_n} $ by a factor of $2.6$ ($-4.2$dB) with respect to NTPv4, but it also obtains smaller confidence intervals and maximum offset values. 

%

\begin{figure}[htp]
\centering
\includegraphics[width=.85\columnwidth]{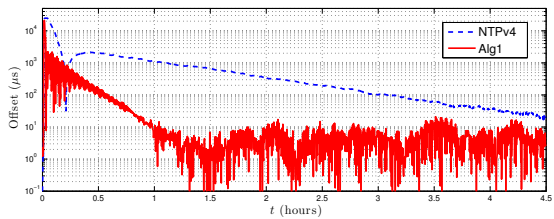}
\caption{Offset values of NTPv4 and Alg1 after a 25ms offset introduced in serv1. }\label{fig:step_ntpvsalg1}
\end{figure}\vspace{-.25cm}

Finally, we investigate the speed of convergence. Starting from both clients synchronized to server serv1, we introduce a 25ms offset. Figure \ref{fig:step_ntpvsalg1} shows how Alg1 is able to converge to a $20\mu s$ range within one hour while NTPv4 needs $4.5$hours to achieve the same synchronization precision.

\vspace{.15cm}
\noindent{\bf Experiment 4 (Comparison with IBM CCT):}
We now proceed to compare the performance of Alg1 with respect to IBM CCT. Notice that unlike IBM CCT, our solution does not perform any previous filtering of the offset samples, the filtering is performed instead by calibrating the parameters which mostly depend on the polling interval $\tau$ chosen. Here we use  $c=0.70$, $\tau=250ms$, $\kp_1=0.1385$, $\kp_2=0.1363$ and $p=0.62$.

\begin{figure}[htp]
       \begin{subfigure}[b]{0.49\columnwidth}
               \centering
               \includegraphics[width=\columnwidth]{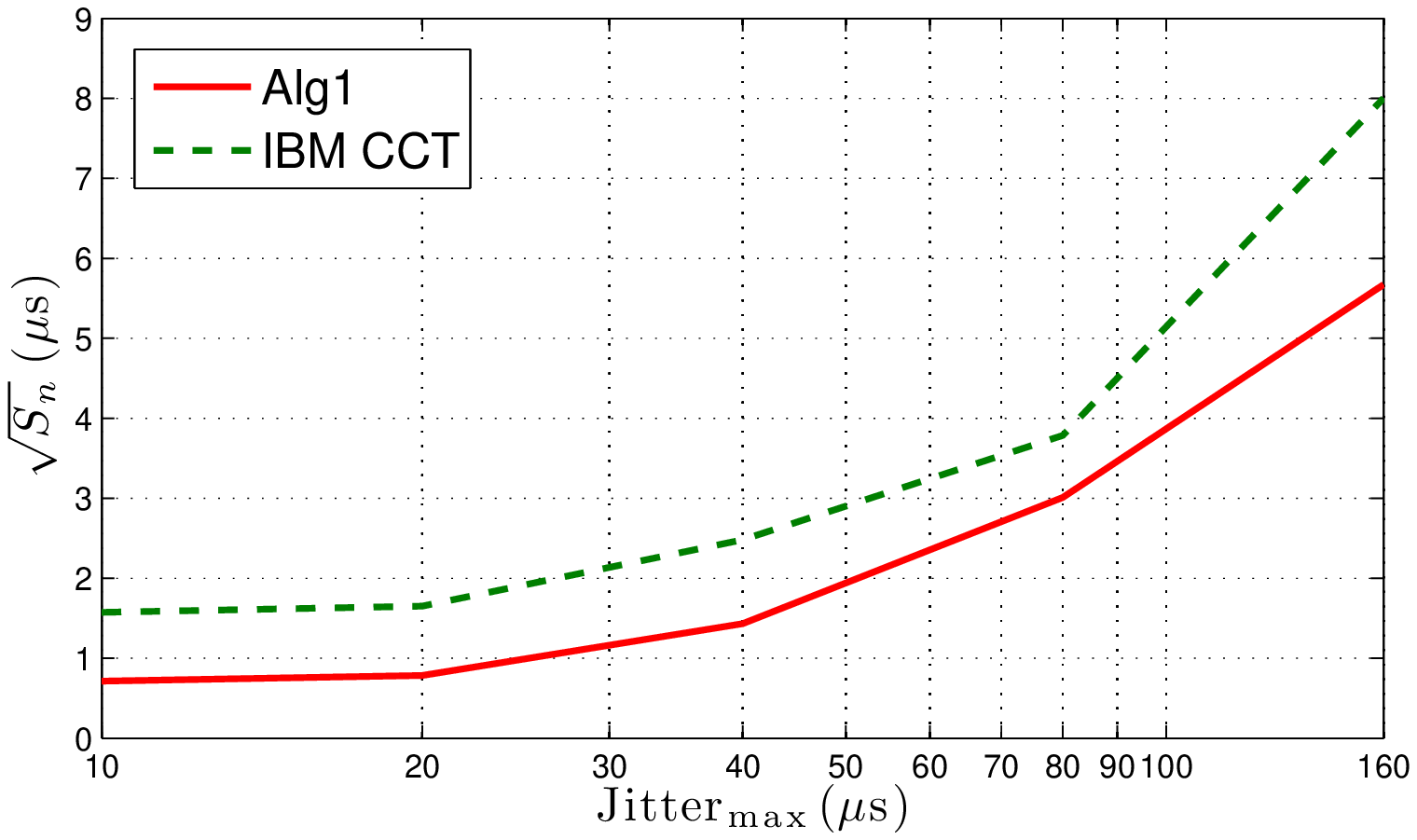}
               \caption{Mean relative deviation $\sqrt{S_n}$}\label{fig:comparison1}
        \end{subfigure}
        \begin{subfigure}[b]{0.49\columnwidth}
               \centering
               \includegraphics[width=\columnwidth]{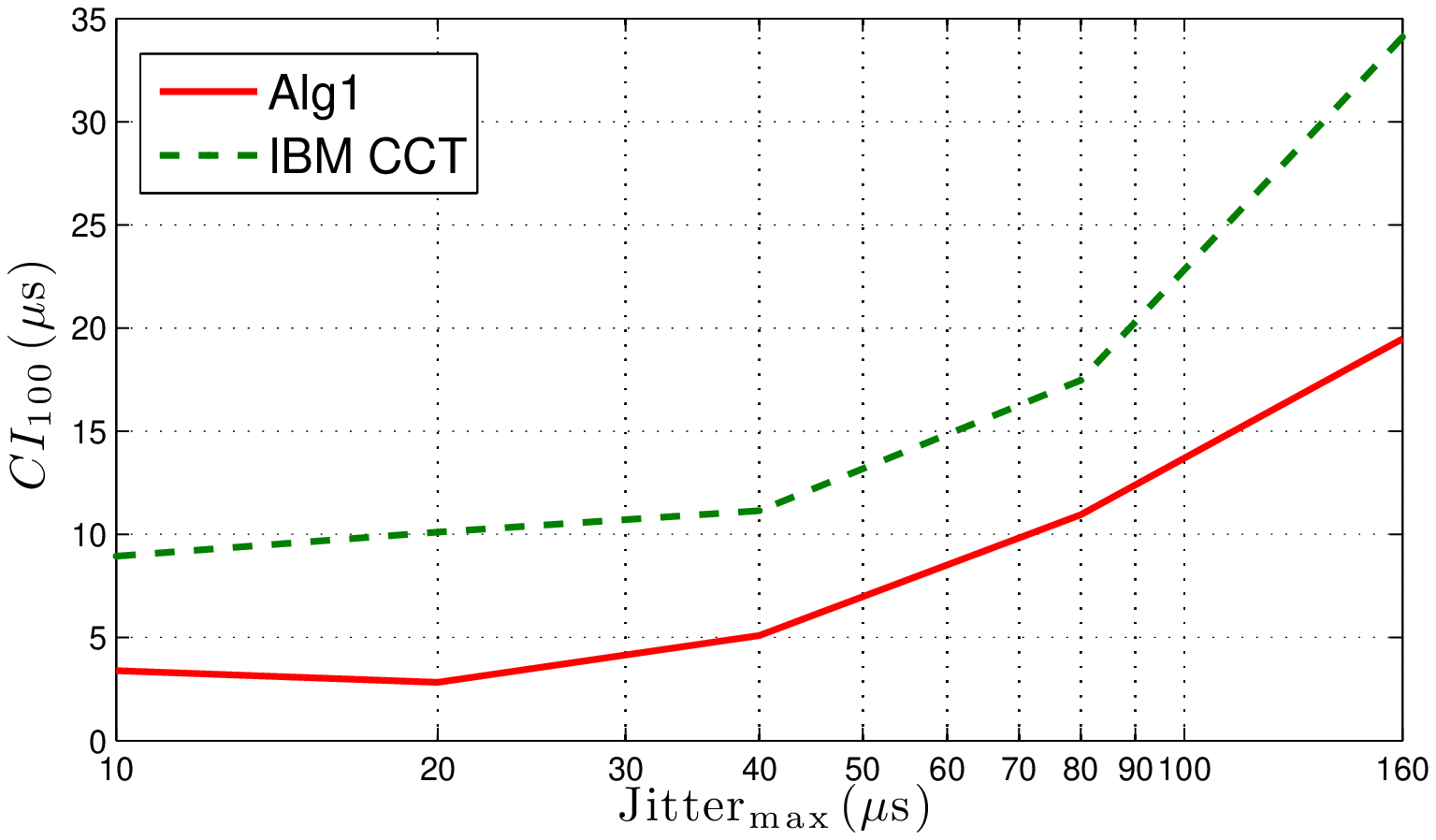}
               \caption{Maximum offset}\label{fig:comparison2}
        \end{subfigure}%
        \vspace{-.25cm}
        \caption{Performance evaluation between our solution (Alg1) and IBM CCT}\label{fig:comparison}
\end{figure}

In Figure \ref{fig:comparison1} we present the mean relative deviation $\sqrt{S_n}$ for two clients connected directly to the leader as the jitter is increased from Jitter$_{\max}=0\mu$s (no jitter) to Jitter$_{\max}=160\mu$s with a granularity of $1\mu$s. The worst case offset is shown in Figure \ref{fig:comparison2}. Each data point is computed using a sample run of 250 seconds.

Our algorithm consistently outperforms IBM CCT in terms of both $\sqrt{S_n}$ and worst case offset. 
The performance improvement is due to two reasons. Firstly, the noise filter used by the IBM CCT algorithm is tailored for noise distributions that are mostly concentrated close to zero with sporadic large errors. However, it does not work properly in cases where the distribution is more homogeneous as in this case. 
Secondly, by choosing $\delta\kp = \kp_1-\kp_2=0.002\ll1$ 
the protocol becomes very robust to offset errors.

\vspace{-.25cm}

\vspace{-.75cm}
\section{Conclusion}\label{sec:conclusions}
This paper presents a clock synchronization protocol that is able to synchronize networked nodes without explicit estimation of the clock skews and steep corrections on the time. Unlike current standards, our protocol is guaranteed to converge even in the presence of timing loops which allows different clients to share timing information and even collectively outperform individual clients. We implemented our solution on a cluster of IBM BladeCenter servers and empirically verified our predictions and our protocol's supremacy over several existing solutions.

%
%

\bibliography{biblio}
\bibliographystyle{IEEEtran}

\appendix


\subsection{ Proof of Lemma \ref{lem:multiplicity} }\label{app:lemma1}
\begin{IEEEproof}
We first compute the characteristic polynomial
\begin{align*}
&\det(\lambda I_{3n}-A) =
\left|
\begin{array}{ccc}
(\lambda-1)I_n  & -\tau R &\0_{n\times n} \\
\kappa_1L & (\lambda-1)I_n & \kappa_2I_n \\
pL & 0 & (\lambda-1+p) I_n
\end{array}
\right| \\
&=(\lambda-1)^n
\left|
\begin{array}{cc}
(\lambda-1)I_n+\frac{\tau \kappa_1}{\lambda-1}LR &  \kappa_2I_n\\
\frac{\tau p}{\lambda-1} LR & (\lambda -1 +p) I_n
\end{array}
\right|\\
&=\det\left( (\lambda-1)^2(\lambda-1+p)I_n +  [(\lambda-1)\kappa_1 \right.  \\
&\left.  +(\kappa_2-\kappa_1)]\tau LR\right) = \prod_{l=1}^{n} g_l(\lambda),
\end{align*}
where $g_l(\lambda)$ is as defined in \eqref{eq:g_l} and we have  iteratively use the determinant property of block matrices
$
\det(A) = \det(A_{11}) \det(A\backslash A_{11})
$
where
$
A=\left[\begin{array}{cc}
A_{11} & A_{12}\\
A_{21} & A_{22}
\end{array}
\right]
$ and $A\backslash A_{11}=A_{22} - A_{21}A_{11}^{-1}A_{12}$ is the Schur complement of $A_{11}$~\cite{horn_matrix_1990}.


Thus, $\lambda=1$ is a double root of the characteristic polynomial if and only if $\kp_1\neq\kp_2$, $p>0$ and $\tau LR$ has a simple zero eigenvalue, i.e. \eqref{eq:nu_condition}.
Now, since $R$ is nonsingular \eqref{eq:nu_condition} must hold for the eigenvalues of $L$ as well, which is in fact true if and only if the directed graph $G(V,E)$ is connected~\cite{xie_consensus_2011}.
\end{IEEEproof}

\subsection{Proof of Lemma \ref{lem:jordan_chains} }\label{app:lemma2}

\begin{IEEEproof}
We start by computing the right Jordan chain. By definition of $\zeta_1$,
$
(A - I)\zeta_1 = 0_n.
$
Thus, if $\zeta_1 = [x^T\;s^T\;y^T]^T$, then the following system of equations must be satisfied
\begin{align}
\tau &R s = \0_n\text{ (a), } -\kappa_1 Lx  - \kappa_2 y =\0_n\text{ (b) }\text{ and }\nonumber\\
 -p &Lx - p y =\0_n\text{ (c).} \label{eq:eq_sys}
\end{align}
Equation (\ref{eq:eq_sys}a) implies $s=0$. Now, since $p>0$, (\ref{eq:eq_sys}c) implies $Lx=-y$, which when substituted in (\ref{eq:eq_sys}b) gives
$
(\kappa_2-\kappa_1)y=\0_n.
$
Thus, since $\kappa_1\neq\kappa_2$, $y=\0_n$ and $x\in\ker(L)$. By choosing $x=\alpha_1\1_n$ (for some $\alpha_1\neq 0$) we obtain 
$
\zeta_1= \alpha_1
\left[
\1_n^T \; \0_n^T \; \0_n^T
\right]^T.
$

Notice that the computation also shows that $\zeta_1$ is the unique eigenvector of $\mu(A)=1$ which implies that there is only one Jordan block of size 2. The second member of the chain, $\zeta_2$, and $\zeta_3$ can be computed similarly by solving
$(A - I_n)\zeta_2 = \zeta_1$ and $(A-(1-p)I_n)\zeta_3=\0_n$. This gives
\begin{align*}
\zeta_2= \left[\begin{array}{c}
\alpha_2\1_n \\ \frac{\alpha_1}{\tau}R^{-1}\1_n \\ \0_n
\end{array}\right] 
\quad\text{ and }\quad 
\zeta_3= \alpha_3\left[\begin{array}{c}
-\frac{\tau\kp_2}{p^2}\1_n \\ \frac{\kp_2}{p}R^{-1}\1_n \\ R^{-1}\1_n
\end{array}\right].
\end{align*} 
In computing $\zeta_3$, we obtain $Lx=0$ and $Rx=-\frac{\tau}{p}s=-\frac{\kp_2\tau}{p^2}y$. $\zeta_3$ follows by taking $y=\alpha_3R^{-1}\1_n$.

The vectors $\eta_1$, $\eta_2$ and $\eta_3$ can be solved in the same way using $\eta_2^T(A - I)=\0_n^T$,
$\eta_1^T(A - I)=\eta_2^T$ and $\eta_3^T(A - (1-p)I)=\0_n^T$. This gives
$\eta_1= \left[
\frac{\beta_2}{\tau}R^{-1}\xi^T \; \beta_1\xi^T \; (-\frac{\kp_2}{p}\beta_1 + \frac{\kp_2}{p^2}\beta_2)\xi^T
\right]^T
$,
$\eta_2= \beta_2\left[
\0_n^T \; \xi^T \; \frac{\kp_2}{p}\xi^T
\right]^T 
\quad\text{ and }\quad
\eta_3= \beta_3\left[
\0_n^T \; \0_n^T \; \xi^T
\right]^T.
$
We set  $\alpha_1=\alpha_2=\alpha_3=1$; this can be done WLOG  provided we still satisfy $\eta_l^T\zeta_l=1$ and $\eta_l^T\zeta_h=0$ for $l\neq h$.
Finally, $\eta_1^T\zeta_1=1$ gives$\beta_2=\gamma\tau$, $\eta_3^T\zeta_3=1$ gives $\beta_3=\gamma$ and $\eta_1^T\zeta_2=0$ gives $\beta_1=-\beta_2=-\gamma\tau$.
\end{IEEEproof}

\subsection{Proof of Theorem \ref{th:convergence} }\label{app:th:convergence}
\begin{IEEEproof}
We first notice that whenever $x(t_k)$ approaches \eqref{eq:synchronization} then
\begin{equation}\label{eq:limit}
\vspace{-.4cm}
\lim_{h\rightarrow\infty} x(t_h) - r^*\1_n (t_h-t_0) = x^*\1_n
\vspace{-.1cm}
\end{equation}

\subsubsection*{Sufficiency}
Since we are under the assumptions of Lemmas \ref{lem:multiplicity} and \ref{lem:jordan_chains} we know that $\mu(A)=1$ has multiplicity $2$ and a Jordan chain of size $2$.
Thus, the Jordan normal form of $A$ is
\begin{equation}\label{eq:jordan_form2}
 A =
\left[ \zeta_1 ... \zeta_{3n}\right]
\left[
\begin{array}{cc}
\begin{array}{crc}
1 & 1 & 0\\
0 & 1 & 0\\
0 & \;\; 0 & 1-p
\end{array}  & \0_{3\times 3(n-1)} \\
\0_{3(n-1)\times 3} & \hat J
\end{array}
\right]
\left[\begin{array}{c} {\eta_1}^T \\\vdots \\ {\eta_{3n}}^T\end{array}\right]
\end{equation}
where $\hat J$ has eigenvalues with spectral radius $\rho(\hat J):=\max_l|\mu_l(\hat J)|<1$. Thus, it follows that
\begin{align}
&\lim_{h\rightarrow \infty} A^h - \zeta_1\eta_1^T - (h\zeta_1 +\zeta_2)\eta_2^T= \lim_{h\rightarrow \infty} \left[ \zeta_1  ... \zeta_{3n}\right]
\label{eq:matrix_limit}\\
& \left[
\begin{array}{cc}
\begin{array}{cc} \0_{2\times 2} &\0_{2\times 1}\\ \0_{1\times 2} & (1-p)^h \end{array} & \0_{2\times(3n-2)} \\
\0_{(3n-2)\times 2} & \hat J^h
\end{array}
\right]
\left[\begin{array}{c} {\eta_1}^T \\\vdots \\ {\eta_{3n}}^T\end{array}\right]
=\0_{3n}\nonumber
\end{align}
where the last equality follows since $(1-p)^h\xrightarrow[h\rightarrow \infty]{}0$ and
$
\Norm{\hat J^h}_\varepsilon\leq \Norm{\hat J}^h_\varepsilon \leq (\rho+\varepsilon)^h\xrightarrow[h\rightarrow \infty]{} 0,
$
where the norm $\Norm{\cdot}_\varepsilon$ is chosen such that $\Norm{A}_\varepsilon=\rho(A) + \varepsilon$~\cite[p. 297, Lemma 5.6.10]{horn_matrix_1990} and $\varepsilon$ is such $\rho(\hat J)+\varepsilon<1$.

Right multiplying \eqref{eq:matrix_limit} with a given initial condition $z_0 = [x_0^T\;s_0^T\;y_0^T]^T$ and using \eqref{eq:zeta} and \eqref{eq:eta} gives
\begin{align}\label{eq:lim2}
\lim_{k\rightarrow \infty}&x_k  -(t_k-t_0)\gamma\1_n\xi^T(s_0 - \frac{\kappa_2}{p}y_0)=\quad\quad\nonumber\\
&=\gamma\1_n\xi^T(R^{-1}x_0+\tau\frac{\kp_2}{p^2}y_0).
\end{align}
Thus, equation \eqref{eq:xs&rs} follows from identifying \eqref{eq:lim2} and \eqref{eq:limit}.

%

\subsubsection*{Necessity} 

The algorithm achieves synchronization whenever  \eqref{eq:limit} holds. Then, it follows from \eqref{eq:system_z} and \eqref{eq:limit} that asymptotically the system behaves according to
\begin{align*}
z_{k} &= \left[
\begin{array}{c}
x_{k} \\ s_{k}\\ y_{k}
\end{array}
\right] =
\left[
\begin{array}{c}
x^*\1_n \\ r^* R^{-1}\1_n\\ \0_n
\end{array}
\right] + k\left[
\begin{array}{c}
\tau r^* \1_n \\ \0_n\\ \0_n 
\end{array}
\right] \\
&=  \left(\tau r^*\zeta_2 + (x^*-\tau r^*)\zeta_1\right) + kr^*\tau \zeta_2.
\end{align*}
Thus, since $P$ is invertible $\zeta_l$ are linearly independent. Therefore, if the system synchronizes for arbitrary initial condition, then it must be the case that the effect of the remaining modes $\mu_l(\Gamma)$ vanishes, which can only happen if for every $\mu_l(\Gamma)\neq 1$, $|\mu_l(\Gamma)|<1$ and the multiplicity of $\mu_l(\Gamma)=1$ is two. Now suppose that either $\kappa_1=\kappa_2$ or $p=0$. Then by Lemma \ref{lem:multiplicity}, the multiplicity of $\mu_l(\Gamma)=1$ is not two which is a contradiction. Thus, we must have $\kappa_1\neq \kappa_2$ and $p>0$ whenever the system synchronizes for arbitrary initial condition.
%
\end{IEEEproof}

\subsection{ Proof of Theorem \ref{th:param_sync} }\label{app:th:param_sync}

\begin{IEEEproof}
We will show that when $G(V,E)$ is connected with $\mu(L)\in\mathds R$, then (i)-(iii) are equivalent to the conditions of Theorem \ref{th:convergence}.

Since, $G(V,E)$ is connected and (i)-(ii) satisfies $p>0$ and $\kp_1\neq\kp_2$, the conditions of  Lemma \ref{lem:multiplicity} are satisfied. Therefore the multiplicity of $\mu(A)=1$ is two and by \eqref{eq:nu_condition} these are the roots of
$g_n(\lambda) = (\lambda-1)^2(\lambda -1 +p),$
 which corresponds to the case $\nu_n=0$.

Thus, to satisfy Theorem \ref{th:convergence} we need to show that the remaining eigenvalues are strictly in the unit circle. This is true for the remaining root of $g_n(\lambda)$ iff (i).


For the remaining $g_l(\lambda)$, this implies that are Schur polynomials. Thus, we will show that $g_l(\lambda)$ is a Schur polynomial if and only if (i)-(iii) hold.
We drop the subindex $l$ for the rest of the proof.

We first transform the Schur stability problem into a Hurwitz stability problem. Consider the change of variable $\lambda = \frac{s+1}{s-1}$. Then $|\lambda|<1$ if and only if $\mathds R[s]<0$.

%

Now, since $\nu>0$ by \eqref{eq:nu_condition}, let
\begin{align*}
P(s) & =\frac{(s-1)^3 }{\delta\kp p\nu}g\left(\frac{s+1}{s-1}\right) =
 s^3 + \left(\frac{2\kappa_1}{\delta\kp p}-3\right)s^2 \\&+ \left( \frac{4}{\delta\kp\nu} + 3 -\frac{4\kappa_1}{\delta\kp p} \right)s + \frac{4(2-p)}{\delta\kp p\nu} +\frac{2\kp_1}{\delta\kp p}-1
\end{align*}
where $\delta\kp = \kp_1 -\kp_2$.

We will apply Hermite-Biehler Theorem  to $P(s)$, but first let us express what $1)$ and $2)$ of Theorem \ref{th:HB} mean here.

Condition $1)$  becomes
\begin{align}\label{eq:HB1}
\frac{2\kappa_1}{\delta\kp p}-3 > 0.
\end{align}

Now let $P^r(\omega)$ and $P^i(\omega)$ be as in  Theorem \ref{th:HB}, i.e. let
\begin{align*}
P^r(\omega)=& -\left(\frac{2\kappa_1}{\delta\kp p}-3\right)\omega^2+\frac{4(2-p)}{\delta\kp p\nu} +\frac{2\kp_1}{\delta\kp p}-1\\
P^i(\omega)=&-\omega^3 + \left( \frac{4}{\delta\kp\nu} + 3 -\frac{4\kappa_1}{\delta\kp p} \right)\omega
\end{align*}

The roots of $P^r(\omega)$ and $P^i(\omega)$ are given by $\omega_0=\pm\sqrt{\omega^r_0}$ and $\omega_0\in\{0,\;\pm\sqrt{\omega^i_0}\}$ respectively, where
\begin{equation}\label{eq:root}
\omega^r_0:=\frac{\frac{4(2-p)}{\delta\kp p\nu} +\frac{2\kp_1}{\delta\kp p}-1}{\frac{2\kappa_1}{\delta\kp p}-3}\text{ and }
\omega^i_0:= \frac{4}{\delta\kp\nu} + 3 -\frac{4\kappa_1}{\delta\kp p}
\end{equation}

Since the roots $P^r(\omega)$ and $P^i(\omega)$ must be real, we must have $\omega_0^r>0$ and $\omega_0^i>0$. Therefore, by monotonicity of the square root, the interlacing condition $2)$ is equivalent to
\begin{equation}\label{eq:HB2}
0<\omega_0^r<\omega_0^i.
\end{equation}
Thus we will show:  (i)-(iii) hold $\iff$ \eqref{eq:HB1} and \eqref{eq:HB2} hold.

It is straightforward to see that using (i) and (ii) we can get \eqref{eq:HB1}. On the other hand,  $\omega_o^i>0$ from \eqref{eq:HB2} together with \eqref{eq:HB1} gives
$
0<\frac{4}{\delta\kp\nu} + 3 -\frac{4\kappa_1}{\delta\kp p} <\frac{4}{\delta\kp\nu}
$, which implies that $\delta\kp>0$, and therefore (ii) follows.

Now using \eqref{eq:HB1} and the definition of $\omega^r_0$ in \eqref{eq:root}, $\omega_0^r>0$ becomes
$
\frac{4(2-p)}{\delta\kp p\nu} +\frac{2\kp_1}{\delta\kp p}-1>0
$
which always holds under (i) and (ii) since the first term is always positive and
$
\frac{2\kp_1}{\delta\kp p}-1>\frac{2\kp_1}{\delta\kp p}-3>0
$
by \eqref{eq:HB1}.

Using \eqref{eq:root}, $\omega_0^r<\omega_0^i$ is equivalent to  
%
\begin{equation}
\nu < \frac{p(\kp_2 - \delta\kp p)}{(\kp_1-\delta\kp p)^2}. \label{eq:inequality}
\end{equation}
Finally, $\nu_l=\mu_l(\tau LR)=\tau \mu_{l}(LR)$. Thus, since \eqref{eq:inequality} should hold $\forall l\in\{1,...,n-1\}$, then
\[
\tau < \min_l \frac{p(\kp_2 - \delta\kp p)}{\mu_l(LR)(\kp_1-\delta\kp p)^2} = \frac{p(\kp_2 - \delta\kp p)}{\mu_{\max}(\kp_1-\delta\kp p)^2}
\]
which is exactly (iii).
\end{IEEEproof}

%
%
%

\end{document}